# POSTERIOR PROPRIETY AND ADMISSIBILITY OF HYPERPRIORS IN NORMAL HIERARCHICAL MODELS[1]

By James O. Berger, William Strawderman and Dejun Tang

*Duke University and SAMSI, Rutgers University and Novartis Pharmaceuticals*

Hierarchical modeling is wonderful and here to stay, but hyperparameter priors are often chosen in a casual fashion. Unfortunately, as the number of hyperparameters grows, the effects of casual choices can multiply, leading to considerably inferior performance. As an extreme, but not uncommon, example use of the wrong hyperparameter priors can even lead to impropriety of the posterior.

For exchangeable hierarchical multivariate normal models, we first determine when a standard class of hierarchical priors results in proper or improper posteriors. We next determine which elements of this class lead to admissible estimators of the mean under quadratic loss; such considerations provide one useful guideline for choice among hierarchical priors. Finally, computational issues with the resulting posterior distributions are addressed.

## 1. Introduction.

1.1. *The model and the problems.* Consider the block multivariate normal situation (sometimes called the "matrix of means problem") specified by the following hierarchical Bayesian model:

$$\mathbf{X} \sim N_p(\boldsymbol{\theta}, \mathbf{I}), \boldsymbol{\theta} \sim N_p(\mathbf{B}, \boldsymbol{\Sigma}_\pi), \tag{1}$$

where

$$\mathbf{X}_{p\times 1} = \begin{pmatrix} \mathbf{X}_1 \\ \mathbf{X}_2 \\ \vdots \\ \mathbf{X}_m \end{pmatrix}, \qquad \boldsymbol{\theta}_{p\times 1} = \begin{pmatrix} \boldsymbol{\theta}_1 \\ \boldsymbol{\theta}_2 \\ \vdots \\ \boldsymbol{\theta}_m \end{pmatrix},$$

Received February 2004; revised July 2004.

[1]Supported by NSF Grants DMS-98-02261 and DMS-01-03265.

$$\mathbf{B}_{p\times 1} = \begin{pmatrix} \boldsymbol{\beta} \\ \boldsymbol{\beta} \\ \vdots \\ \boldsymbol{\beta} \end{pmatrix}, \qquad \boldsymbol{\Sigma}_{\pi p\times p} = \begin{pmatrix} \mathbf{V} & \mathbf{0} & \cdots & \mathbf{0} \\ \mathbf{0} & \mathbf{V} & \cdots & \mathbf{0} \\ \vdots & \vdots & \ddots & \vdots \\ \mathbf{0} & \mathbf{0} & \cdots & \mathbf{V} \end{pmatrix},$$

where the $\mathbf{X}_i$ are $k \times 1$ observation vectors, $k \geq 2$, the $\boldsymbol{\theta}_i$ are $k \times 1$ unknown mean vectors, $\boldsymbol{\beta}$ is a $k \times 1$ unknown "hyper-mean" vector and $\mathbf{V}$ is an unknown $p \times p$ "hyper-covariance matrix." This is more commonly written as, for $i = 1, 2, \ldots, m$ and independently, $\mathbf{X}_i \sim N_k(\boldsymbol{\theta}_i, \mathbf{I})$, $\boldsymbol{\theta}_i \sim N_k(\boldsymbol{\beta}, \mathbf{V})$. Note that $p = mk$. Efron and Morris [16, 17] introduced the study of this model from an empirical Bayes perspective. Today, it is more common to analyze the model from a hierarchical Bayesian perspective (cf. [2, 18]), based on choice of a hyperprior $\pi(\boldsymbol{\beta}, \mathbf{V})$. Such hyperpriors are often chosen quite casually, for example, constant priors or the "nonhierarchical independence Jeffreys prior" (see Section 1.2) $|\mathbf{V}|^{-(k+1)/2}$. In this paper we formally study properties of such choices.

The first issue that arises when using improper hyperpriors is that of propriety of the resulting posterior distributions (cf. [38]). In Section 2 we discuss choices of $\pi(\boldsymbol{\beta}, \mathbf{V})$ which yield proper posterior distributions. That this is of importance is illustrated by the fact that we have seen many instances of use of $|\mathbf{V}|^{-(k+1)/2}$ for similar situations, even though it is known to generally yield an improper posterior distribution when used as a hyperprior (see Section 2).

A more refined question, from the decision-theoretic point of view, is that of choosing hyperpriors so that the resulting Bayes estimators, for a specified loss function, are admissible. The particular version of this problem that we will study is that of estimating $\boldsymbol{\theta}$ by its posterior mean $\boldsymbol{\delta}^{\pi}(\mathbf{x})$, under the quadratic loss

$$L(\boldsymbol{\theta}, \boldsymbol{\delta}^{\pi}) = (\boldsymbol{\theta} - \boldsymbol{\delta}^{\pi})^t \mathbf{Q} (\boldsymbol{\theta} - \boldsymbol{\delta}^{\pi}), \tag{2}$$

where $\mathbf{Q}$ is a known positive-definite matrix. The performance of an estimator $\boldsymbol{\delta}$ will be evaluated by the usual frequentist risk function

$$R(\boldsymbol{\theta}, \boldsymbol{\delta}) = E_{\boldsymbol{\theta}}^{\mathbf{X}}[L(\boldsymbol{\theta}, \boldsymbol{\delta}(\mathbf{X}))]. \tag{3}$$

The estimator $\boldsymbol{\delta}$ is *inadmissible* if there exists another estimator with risk function nowhere bigger and somewhere smaller. If no such better estimator exists, $\boldsymbol{\delta}$ is *admissible*.

In Section 3 conditions on $\pi(\boldsymbol{\beta}, \mathbf{V})$ are presented under which the Bayes estimator $\boldsymbol{\delta}^{\pi}$ is admissible and inadmissible. The motivation for looking at this problem is not that this specific decision-theoretic formulation is necessarily of major practical importance. The motivation is, instead, that use of "objective" improper priors in hierarchical modeling is of enormous practical importance, yet little is known about which such priors are good or bad. The most successful approach to evaluation of objective improper priors has been to study the frequentist properties of the ensuing Bayes procedures (see [3] for discussion and many references). In particular, it is important that the prior distribution not be too diffuse, and study of admissibility is the most powerful tool known for detecting an over-diffuse prior. Also see [10] for general discussion of the utility of the decision-theoretic perspective in modern statistical inference.

The results in the paper generalize immediately to the case where the identity covariance matrix $\mathbf{I}$ for the $\mathbf{X}_i$ is replaced by a known positive-definite covariance matrix, but for notational simplicity we only consider the identity case. More generally, the motivation for this study is to obtain insight into the choice of hyperpriors in multivariate hierarchical situations. The possibilities for normal hierarchical modeling are endless, and it is barely conceivable that formal results about posterior propriety and admissibility can be obtained in general. The hope behind this study is that what is learned in this specific multivariate hierarchical model can provide guidance in more complex hierarchical models.

1.2. *The hyperprior distributions being studied.* We will study hyperprior densities of the form

$$\pi(\boldsymbol{\beta}, \mathbf{V}) = \pi(\boldsymbol{\beta})\pi(\mathbf{V}).$$

For $\mathbf{V}$, we will study priors that satisfy the following condition, where $d_1 > d_2 > \cdots > d_k > 0$ are the eigenvalues of $\mathbf{V}$.

CONDITION 1. For $0 \le l \le 1$,

$$\frac{C_1}{|\mathbf{I}+\mathbf{V}|^{(a_2-a_1)}|\mathbf{V}|^{a_1}[\prod_{i<j}(d_i-d_j)]^{(1-l)}} \le \pi(\mathbf{V}) \le \frac{C_2}{|\mathbf{I}+\mathbf{V}|^{(a_2-a_1)}|\mathbf{V}|^{a_1}[\prod_{i<j}(d_i-d_j)]^{(1-l)}},$$

where $C_1$ and $C_2$ are positive constants and $|\mathbf{A}|$ denotes the determinant of $\mathbf{A}$. Many common noninformative priors satisfy this condition, including:

*Constant prior.* $\pi(\mathbf{V}) = 1$; here $a_1 = a_2 = 0$ and $l = 1$.

*Nonhierarchical independence Jeffreys prior.* $\pi(\mathbf{V}) = |\mathbf{V}|^{-(k+1)/2}$; here $a_1 = a_2 = (k+1)/2$ and $l = 1$.

*Hierarchical independence Jeffreys prior.* $\pi(\mathbf{V}) = |\mathbf{I}+\mathbf{V}|^{-(k+1)/2}$; here $a_1 = 0, a_2 = (k+1)/2$ and $l = 1$.

*Nonhierarchical reference prior.* $\pi(\mathbf{V}) = [|\mathbf{V}| \prod_{i<j}(d_i - d_j)]^{-1}$; here $a_1 = a_2 = 1$ and $l = 0$. (See [40].)

*Hierarchical reference priors.*

(a) $\pi(\mathbf{V}) = [|\mathbf{I} + \mathbf{V}| \prod_{i<j}(d_i - d_j)]^{-1}$; here $a_1 = 0$, $a_2 = 1$ and $l = 0$.

(b) $\pi(\mathbf{V}) = [|\mathbf{V}|^{-(2k-1)/(2k)} \prod_{i<j}(d_i - d_j)]^{-1}$; here $a_1 = a_2 = (2k-1)/(2k)$ and $l = 0$.

We have already alluded to the nonhierarchical independence Jeffreys prior, which formally is the Jeffreys prior for a covariance matrix in a nonhierarchical setting with given mean. Unfortunately, this prior seems to be commonly used for covariance matrices at any level of a hierarchy, typically yielding improper posteriors, as will be seen in Section 2. Those who recognize the problem often instead use the constant prior, or the hierarchical independence Jeffreys prior, which arises from considering the "marginal model" formed by integrating over $\boldsymbol{\beta}$ in the original model and computing the independence Jeffreys prior for this marginal model.

Similarly, the nonhierarchical reference prior yields an improper posterior in hierarchical settings (shown in Section 2). The two versions of hierarchical reference priors given above arise from quite different perspectives. Prior (a) arises from considering the marginal model formed by integrating over $\boldsymbol{\beta}$ in the original model, and applying the Yang and Berger [40] reference prior formula to the covariance matrix $\mathbf{I} + \mathbf{V}$ that arises in the marginal model. (The differences of eigenvalues for this matrix are the same as the differences of the eigenvalues for $\mathbf{V}$.) Prior (b) arises from a combination of computational and admissibility considerations that are summarized in Sections 1.5 and 1.6, respectively.

Note that if the covariance matrix for the $\mathbf{X}_i$ were a known $\boldsymbol{\Sigma}$, instead of $\mathbf{I}$, then $\mathbf{I}$ in the above priors would be replaced by $\boldsymbol{\Sigma}$. It could not then be said, however, that the reference prior formula is that which would result from applying the Yang and Berger [40] reference prior formula to the covariance matrix $\boldsymbol{\Sigma} + \mathbf{V}$ that arises in the marginal model, since the differences of eigenvalues of this matrix will no longer equal the differences of the eigenvalues of $\mathbf{V}$.

Three commonly considered priors for the hyperparameter $\boldsymbol{\beta}$ are:

*Case* 1. *Constant prior.* $\pi(\boldsymbol{\beta}) = 1$.

*Case* 2. *Conjugate prior.* $\pi(\boldsymbol{\beta})$ is $N_k(\boldsymbol{\beta}^0, \mathbf{A})$, where $\boldsymbol{\beta}^0$ and $\mathbf{A}$ are subjectively specified.

*Case* 3. *Hierarchical prior.* $\pi(\boldsymbol{\beta})$ is itself given in two stages:

$$\boldsymbol{\beta}|\lambda \sim N_k(\boldsymbol{\beta}^0, \lambda\mathbf{A}), \qquad \lambda \sim \pi(\lambda), \lambda > 0, \tag{4}$$

where $\boldsymbol{\beta}^0$ and $\mathbf{A}$ are again specified, and $\pi(\lambda)$ satisfies:

CONDITION 2.

(i) $\int_0^c \pi(\lambda)\,d\lambda < \infty$ for $c > 0$;
(ii) $\pi(\lambda) \sim C\lambda^{-b}$ $(b \geq 0)$ as $\lambda \to \infty$ for some constant $C > 0$.

As discussed in [7], an important example of a Case 3 prior is obtained by choosing

$$\pi(\lambda) \propto \lambda^{-b} e^{-c/\lambda},$$

that is, an inverse Gamma$(b-1, c^{-1})$ density. This clearly satisfies Condition 2, and the resulting prior for $\boldsymbol{\beta}$ is

$$\pi(\boldsymbol{\beta}) = \int \pi(\boldsymbol{\beta}|\lambda)\pi(\lambda)\,d\lambda \propto \left[1 + \frac{1}{2c}(\boldsymbol{\beta} - \boldsymbol{\beta}^0)^t \mathbf{A}^{-1}(\boldsymbol{\beta} - \boldsymbol{\beta}^0)\right]^{-(k/2+b-1)},$$

which is a multivariate $t$-distribution with median $\boldsymbol{\beta}^0$, scale matrix proportional to $\mathbf{A}$ and $2(b-1)$ degrees of freedom. We will be particularly interested in the improper version of this prior with $c = 1/2$, $\boldsymbol{\beta}^0 = \mathbf{0}$, $\mathbf{A} = \mathbf{I}$ and $b = 1/2$, corresponding to

$$\pi(\boldsymbol{\beta}) \propto [1 + \|\boldsymbol{\beta}\|^2]^{-(k-1)/2}. \tag{5}$$

1.3. *Related literature.* Hierarchical Bayesian analysis has been widely applied to many theoretical and practical problems (cf. [8, 11, 21, 23]). Results and many references to decision-theoretic analysis of hierarchical Bayesian models can be found in [5, 6, 7, 18, 31]. Reference [7] considered the following hierarchical normal model: $\mathbf{X} = (X_1, X_2, \ldots, X_p)^t \sim N_p(\boldsymbol{\theta}, \boldsymbol{\Sigma})$, with $\boldsymbol{\Sigma}$ being a known positive-definite matrix. The paper considered the common two-stage prior distribution for $\boldsymbol{\theta}$ given by $\boldsymbol{\theta} \sim N_p(\beta\underline{1}, \sigma_\pi^2\mathbf{I}), (\beta, \sigma_\pi^2) \sim \pi_1(\sigma_\pi^2)\pi_2(\beta)$, where $\underline{1}$ is the $p$-vector of 1's and $\mathbf{I}$ is the identity matrix, and presented choices of $\pi_1(\sigma_\pi^2)$ and $\pi_2(\beta)$ which yield proper posteriors and admissible Bayes estimators under quadratic loss. This is thus the special case of our model where $k = 1$, and this paper can be viewed as an extension of those results to the vector mean problem (and, hence, our restriction to $k \geq 2$).

The more general decision-theoretic background of this paper is the huge literature on shrinkage estimation, initiated by the demonstration in [33] that the usual estimator for the mean of a multivariate normal distribution is not admissible when $p \geq 3$. This huge literature can be accessed from, for instance, [36].

The key to the admissibility and inadmissibility results presented in this paper is the fundamental paper [9], which provided the crucial insight to allow determination of admissibility and inadmissibility of Bayes estimators.

1.4. *A transformation and revealed concern.* It is convenient, for both intuitive and technical reasons, to write $\mathbf{V} = \mathbf{H}^t\mathbf{D}\mathbf{H}$, where $\mathbf{H}$ is the matrix of eigenvectors corresponding to $\mathbf{D} = \mathrm{diag}(d_1, d_2, \ldots, d_k)$, such that $\mathbf{H}^t\mathbf{H} = \mathbf{I}$. Indeed, we will make the change of variables from $\mathbf{V}$ to $(\mathbf{D}, \mathbf{H})$, and rewrite the prior as

$$\pi(\mathbf{V})\, d\mathbf{V} = \pi(\mathbf{H}, \mathbf{D}) I_{[d_1 > d_2 > \cdots > d_k]}\, d\mathbf{D}\, d\mathbf{H};$$

here $d\mathbf{V} = \prod_{i \le j} dV_{ij}$, $d\mathbf{D} = \prod_{i=1}^k dd_i$, $d\mathbf{H}$ denotes the invariant Haar measure over the space of orthonormal matrices and $I_{[d_1 > d_2 > \cdots > d_k]}$ denotes the indicator function over the specified set. (Because of Condition 1, equality of any eigenvalues has measure 0.)

From [19], the functional relationship between $\pi(\mathbf{V})$ and $\pi(\mathbf{H}, \mathbf{D})$ is

$$\pi(\mathbf{H}, \mathbf{D}) = \pi(\mathbf{H}^t\mathbf{D}\mathbf{H}) \prod_{i<j} (d_i - d_j).$$

Thus Condition 1 becomes

CONDITION 1′. For $0 \le l \le 1$,

$$\frac{C_1[\prod_{i<j}(d_i - d_j)]^l}{|\mathbf{I} + \mathbf{D}|^{(a_2 - a_1)} |\mathbf{D}|^{a_1}} \le \pi(\mathbf{H}, \mathbf{D}) \le \frac{C_2[\prod_{i<j}(d_i - d_j)]^l}{|\mathbf{I} + \mathbf{D}|^{(a_2 - a_1)} |\mathbf{D}|^{a_1}}.$$

Under this transformation, the common objective priors for $\mathbf{V}$ are as follows:

1. The constant prior is now $\pi(\mathbf{H}, \mathbf{D}) = \prod_{i<j}(d_i - d_j)$.
2. The nonhierarchical independence Jeffreys prior is $\pi(\mathbf{H}, \mathbf{D}) = |\mathbf{D}|^{-(k+1)/2} \times \prod_{i<j}(d_i - d_j)$.
3. The hierarchical independence Jeffreys prior is $\pi(\mathbf{H}, \mathbf{D}) = |\mathbf{I} + \mathbf{D}|^{-(k+1)/2} \times \prod_{i<j}(d_i - d_j)$.
4. The nonhierarchical reference prior is $\pi(\mathbf{H}, \mathbf{D}) = |\mathbf{D}|^{-1}$.
5. The hierarchical reference priors are (a) $\pi(\mathbf{H}, \mathbf{D}) = |\mathbf{I} + \mathbf{D}|^{-1}$ and (b) $\pi(\mathbf{H}, \mathbf{D}) = |\mathbf{D}|^{-(2k-1)/(2k)}$.

This transformation reveals a significant difficulty of any prior that can be written as a function of $|\mathbf{V}|$: in the $(\mathbf{H}, \mathbf{D})$ space, such priors contain the factor $\prod_{i<j}(d_i - d_j)$, which gives low mass to close eigenvalues, and hence effectively forces the eigenvalues apart. (The effective prior on $\mathbf{H}$ is just constant, which is natural since $\mathbf{H}$ ranges over a compact space, and hence has no effect on the eigenvalues.) This is contrary to common intuition, in that one is often debating between choice of a covariance matrix with equal eigenvalues or choice of an arbitrary covariance matrix; if anything this would suggest that one should choose a prior that pushes the eigenvalues closer together.

This intuition also receives support from the frequentist literature. The independence Jeffreys prior (and often-employed modifications such as $|\mathbf{I}+\mathbf{V}|$) are of this suspicious form and, when used at the first level of a normal model, result in estimates of $\mathbf{V}$ that are proportional to $\mathbf{S}$, the sample covariance matrix. The frequentist literature, starting with [34] and continuing with such works as [22, 26, 27, 30, 40], shows that $\mathbf{S}$ has eigenvalues that are too disperse and that shrinking the eigenvalues of $\mathbf{S}$ together is necessary for good performance. Since multiples of $\mathbf{S}$ arise as Bayes estimators for priors of the "suspicious" form, there appears to be a direct analogy between what frequentists observed about $\mathbf{S}$ and the concern that these priors force the eigenvalues apart.

In contrast to this behavior, the reference and hierarchical reference priors do not contain the term $\prod_{i<j}(d_i - d_j)$ in the transformed space, and hence are neutral with respect to expansion or shrinkage of the eigenvalues. Interestingly, in [40] (see also [37]), it is shown that the Bayes estimators arising from the reference prior (in the nonhierarchical model) behave very similarly to the Stein [34] and Haff [22] estimators, suggesting that such neutral behavior is natural for frequentist estimators—that is, that shrinking the eigenvalues of $\mathbf{S}$ corresponds to a Bayesian prior that is neutral about the eigenvalues. (It should be noted that, in the more recent Bayesian literature, aggressive shrinkage of eigenvalues, correlations or other features of the covariance matrix is entertained; cf. [14, 15, 25] and the references therein. This may well be desirable in many practical situations, but is more aggressive in its prior assumptions than the objective priors we consider.)

1.5. *Computation.* Hierarchical models are typically handled today by Gibbs sampling, possibly with rejection or Metropolis–Hastings steps in the Gibbs sampler (cf. [12, 32]). We briefly indicate considerations in utilizing the priors discussed in Section 1.2 within such computational frameworks.

Use of the Case 1 (constant) or Case 2 (normal) priors for $\boldsymbol{\beta}$ causes no difficulties; sampling of $\boldsymbol{\beta}$ can simply be carried out with a Gibbs step, as its full conditional will be a normal distribution. The Case 3 prior is almost as easy to utilize, because of its representation as a mixture of normals. Indeed, one purposely introduces the latent variable $\lambda$ having the density in (4); sampling of $\boldsymbol{\beta}$ is then done from its full conditional—also given $\lambda$—which is normal, with $\lambda$ then being sampled from its full conditional

$$\pi(\lambda|\boldsymbol{\beta}) \propto \frac{1}{\lambda^{(b+k/2)}} \exp\left(-\frac{1}{\lambda}\left[c + \frac{1}{2}(\boldsymbol{\beta}-\boldsymbol{\beta}^0)^t \mathbf{A}^{-1}(\boldsymbol{\beta}-\boldsymbol{\beta}^0)\right]\right),$$

that is, an inverse Gamma$(b-1+k/2, [c+\frac{1}{2}(\boldsymbol{\beta}-\boldsymbol{\beta}^0)^t\mathbf{A}^{-1}(\boldsymbol{\beta}-\boldsymbol{\beta}^0)]^{-1})$ density. In particular, the recommended default hyperprior $\pi(\boldsymbol{\beta}) \propto [1+\|\boldsymbol{\beta}\|^2]^{-(k-1)/2}$ is handled as above, sampling $\lambda$ from the inverse Gamma$((k-1)/2, 2/[1+\|\boldsymbol{\beta}\|^2])$ density.

Dealing with the hyper-covariance matrix $\mathbf{V}$ is not as easy (cf. [13]), except for the constant prior $\pi(\mathbf{V}) = 1$, for which the full conditional of $V$ is simply an inverse Wishart distribution; alas, this is not a desirable prior in other respects. More attractive is the hierarchical independence Jeffreys prior $\pi(\mathbf{V}) = |\mathbf{I}+\mathbf{V}|^{-(k+1)/2}$. Defining $\mathbf{W}(\boldsymbol{\theta},\boldsymbol{\beta}) = \sum_{i=1}^m(\boldsymbol{\theta}_i-\boldsymbol{\beta})(\boldsymbol{\theta}_i-\boldsymbol{\beta})^t$, the resulting full conditional for $\mathbf{V}$ can be written

$$\pi(\mathbf{V}|\boldsymbol{\theta},\boldsymbol{\beta}) \propto \frac{1}{|\mathbf{I}+\mathbf{V}|^{(k+1)/2}|\mathbf{V}|^{m/2}} \exp\left(-\frac{1}{2}\operatorname{tr}(\mathbf{V}^{-1}\mathbf{W}(\boldsymbol{\theta},\boldsymbol{\beta}))\right),$$

which, unfortunately, is not of closed form. Still, one can easily sample from this full conditional using the following accept-reject sampling algorithm:

*Propose* a candidate $\mathbf{V}^*$ from the inverse Wishart $(m, \mathbf{W}(\boldsymbol{\theta},\boldsymbol{\beta}))$ density

$$g(\mathbf{V}|\mathbf{S}) \propto \frac{1}{|\mathbf{V}|^{m/2+(k+1)/2}} \cdot \exp\left(-\frac{1}{2}\operatorname{tr}(\mathbf{V}^{-1}\mathbf{W}(\boldsymbol{\theta},\boldsymbol{\beta}))\right). \tag{6}$$

*Accept* the candidate with probability $P = (|\mathbf{V}|/|\mathbf{I}+\mathbf{V}|)^{(k+1)/2}$, returning to the proposal step if the candidate is rejected, and moving on to another full conditional if it is accepted.

For large $V$ or $m$ or small dimension $k$, the acceptance probability will be quite high.

When using the hierarchical independence Jeffreys prior, one can gain efficiency by working with the marginal distribution of $\mathbf{V}$, instead of the full conditionals. This is particularly convenient in the Case 1 scenario, where the overall posterior distribution can be written $\pi(\boldsymbol{\theta}|\mathbf{V},\mathbf{x})\pi(\mathbf{V}|\mathbf{x})$, the first posterior being a normal distribution, and hence trivial to sample from, and the marginal posterior of $\mathbf{V}$ being proportional to the integrand in the first expression of Lemma 2.1, namely

$$\pi(\mathbf{V}|\mathbf{x}) \propto \frac{1}{|\mathbf{I}+\mathbf{V}|^{(m+k)/2}} \exp\left(-\frac{1}{2}\sum_{i=1}^m(\mathbf{x}_i-\bar{\mathbf{x}})^t(\mathbf{I}+\mathbf{V})^{-1}(\mathbf{x}_i-\bar{\mathbf{x}})\right).$$

As discussed in [18] (although they utilized the constant prior for $\mathbf{V}$), one can construct a rejection sampler for $\mathbf{V}$ by simply generating $\mathbf{B} = (\mathbf{I}+\mathbf{V}^*)$ from the inverse Wishart $(m+k, \sum_{i=1}^m(\mathbf{x}_i-\bar{\mathbf{x}})(\mathbf{x}_i-\bar{\mathbf{x}})^t)$ density, accepting the candidate $\mathbf{V}^* = \mathbf{B}-\mathbf{I}$ if it is positive definite and returning to generate a new $\mathbf{B}$ if it is not. This will have a reasonable acceptance probability if $\mathbf{V}$ is large or $m$ is large.

For the hierarchical reference priors, it seems that Metropolis–Hastings must be used to sample from the full conditionals. The "standard" approach is that utilized in [40] and [28]. In this approach one first performs the exponential matrix transform of $\mathbf{V}$, which translates the set of positive-definite matrices into unconstrained Euclidean space. Then a hit-and-run Metropolis–Hastings algorithm is employed to produce the Markov chain. This algorithm can be directly utilized here, requiring only the change in the acceptance probability induced by using the hierarchical reference priors instead of the nonhierarchical reference priors.

Since a Metropolis–Hastings step is required anyway for the hierarchical reference priors, one can again gain efficiency by working with the marginal distributions of $\mathbf{V}$, instead of the full conditionals. Taking the Case 3 situation for illustration, one uses the posterior form

$$\pi(\boldsymbol{\theta}|\boldsymbol{\beta},\mathbf{V},\mathbf{x})\pi(\boldsymbol{\beta}|\mathbf{V},\lambda,\mathbf{x})\pi(\mathbf{V},\lambda|\mathbf{x}),$$

where the first two posteriors are simply normal distributions, and hence trivial to sample, and the marginal posterior of $(\mathbf{V},\lambda)$ is proportional to the integrand in the first expression of Lemma 2.3, that is,

$$\begin{aligned}\pi(\mathbf{V},\lambda|\mathbf{x}) \propto {} & \frac{1}{|\mathbf{I}+\mathbf{V}|^{(m-1)/2}|\mathbf{I}+\mathbf{V}+m\lambda\mathbf{A}|^{1/2}} \\ & \times \exp\left(-\frac{1}{2}\sum_{i=1}^{m}(\mathbf{x}_i-\bar{\mathbf{x}})^t(\mathbf{I}+\mathbf{V})^{-1}(\mathbf{x}_i-\bar{\mathbf{x}})\right) \\ & \times \exp\left(-\frac{1}{2}m\bar{\mathbf{x}}^t(\mathbf{I}+\mathbf{V}+m\lambda\mathbf{A})^{-1}\bar{\mathbf{x}}\right)\pi(\mathbf{V})\pi(\lambda).\end{aligned}$$

One proceeds by applying the exponential matrix transform to $\mathbf{V}$ and then running a hit-and-run algorithm for the transformed $\mathbf{V}$ and $\lambda$. For each $(\mathbf{V},\lambda)$ in the chain (or probably better—for, say, every 100th in the chain) one can then generate $\boldsymbol{\beta}$ from the normal $\pi(\boldsymbol{\beta}|\mathbf{V},\lambda,\mathbf{x})$ and then $\boldsymbol{\theta}$ from the normal $\pi(\boldsymbol{\theta}|\boldsymbol{\beta},\mathbf{V},\mathbf{x})$.

If one wishes to stick to Gibbs sampling for the hierarchical reference priors (as would be the case, e.g., if one were working with a complex model for which marginalization could not be carried out), and further desires an easy-to-code algorithm, one could use Metropolis–Hastings on the full conditional for $\mathbf{V}$ with the proposal in (6). (For justification as to why this is the best inverse Wishart proposal, see [39].) The acceptance probabilities for the (a) and (b) versions of the hierarchical reference prior would then be, respectively,

$$\min\left\{1,\frac{\prod_{i<j}(d_i^*-d_j^*)}{\prod_{i<j}(d_i-d_j)}\cdot\frac{|\mathbf{I}+\mathbf{V}^*|\,|\mathbf{V}|^{(k+1)/2}}{|\mathbf{I}+\mathbf{V}|\,|\mathbf{V}^*|^{(k+1)/2}}\right\},$$

$$\min\left\{1,\frac{\prod_{i<j}(d_i^*-d_j^*)}{\prod_{i<j}(d_i-d_j)}\cdot\frac{|\mathbf{V}|^{(k-1+k^{-1})/2}}{|\mathbf{V}^*|^{(k-1+k^{-1})/2}}\right\}.$$

Note, also, that it is generally best to iterate a number of times on this Metropolis step, keeping only the last value, before moving on to another full conditional (as this step is considerably less efficient than the others).

For small $k$ or large $m$, this simple approach will work reasonably well. For instance, in a simulation reported in detail in [39], the average number of Metropolis iterations before a move occurred was as indicated in Table 1. Since the proposal moves widely over the parameter space, a Metropolis scheme that moves at least once in every 10 iterations is often acceptable; thus, for $m \leq 30$, one can use this algorithm to do the calculation with $k$ up to 7. With larger $m$, such as 100, the algorithm is still acceptable for $k = 15$. When this scheme is not efficient enough, the exponential matrix transform hit-and-run approach mentioned above has proven to be very effective (but harder to program).

1.6. *Summary and generalizations.* The results in the paper require significant technical machinery. This machinery is not necessary for understanding the basic conclusions, so we present the most important conclusions and potential generalizations here. Note that the conclusions depend on intuitive appeal (e.g., Section 1.4), posterior propriety (Section 2), admissibility (Section 3) and computational simplicity (Section 1.5).

None of the priors on $\boldsymbol{\beta}$ significantly affects posterior propriety, or caused difficulties in the posterior computation. Hence admissibility is the most important criterion for deciding between them. It seems that use of the constant prior $\pi(\boldsymbol{\beta}) = 1$ results in inadmissibility, except for the case $k = 2$. (This is, of course, not a surprise, in that two dimensions is typically the cut-off for admissibility with constant priors on means.) The Case 2 conjugate prior is, perhaps, reasonable, if one has subjective information about $\boldsymbol{\beta}$. Among the Case 3 default priors, the prior $\pi(\boldsymbol{\beta}) \propto [1 + \|\boldsymbol{\beta}\|^2]^{-(k-1)/2}$ is excellent from the perspective of admissibility for all $k$, and is the prior that we actually recommend for default use. Part of the motivation here is the many studies that have shown the great success of these mixture-of-normals priors in shrinkage estimation in particular (cf. [1, 20, 35]), and robust Bayesian estimation in general (cf. [4]). There is the caveat, however, that this prior should probably only be applied when the $\beta_i$ are roughly "exchangeable," which might well require some reparameterization to ensure. Note that we even recommend use of this prior when $k = 2$. It is often thought that shrinkage should only be used when $k \geq 3$, but it can be used to practical advantage even when $k = 2$ (even though there are no longer uniform dominance results).

TABLE 1
*Average number of nonmoves*

| | ***m*** | | | |
|---|---|---|---|---|
| ***k*** | **20** | **30** | **50** | **100** |
| 3 | 6.89 | 4.92 | 2.14 | 1.06 |
| 5 | 9.83 | 5.74 | 2.96 | 1.21 |
| 7 | 13.52 | 8.50 | 4.03 | 2.27 |
| 10 | 18.74 | 10.86 | 5.42 | 3.46 |
| 12 | 33.67 | 19.63 | 7.61 | 5.07 |
| 15 | 127.35 | 42.98 | 17.89 | 9.36 |

Of considerably more importance than the prior on $\boldsymbol{\beta}$ is the prior on $\mathbf{V}$. The two priors for $\mathbf{V}$ that we have seen most commonly used in practice are the constant prior (or, equivalently, a "vague proper inverse Wishart" prior) and the nonhierarchical Jeffreys prior (or a vague proper inverse Wishart equivalent). Use of the nonhierarchical Jeffreys prior is simply a mistake, in that it results in an improper posterior (and use of the vague proper inverse Wishart equivalent is no better, in that it essentially yields a posterior with almost all its mass in a spike near $\mathbf{V} = \mathbf{0}$). The constant prior requires $m$, the number of blocks, to be about $2k$ in order to achieve posterior propriety. Intuitively, at most $k$ blocks are needed for identifiability of $\mathbf{V}$, so this is a strong indication of the inadequacy of the constant prior. In this regard, the hierarchical independence Jeffreys prior $\pi(\mathbf{V}) = |\mathbf{I} + \mathbf{V}|^{-(k+1)/2}$ requires only $k$ blocks ($k + 1$ if the constant prior on $\boldsymbol{\beta}$ is used) for posterior propriety.

We were not able to establish any admissibility results for these priors, but Tatsuya Kubokawa (private communication) has been able to show by different techniques that the $l = 0$ prior results in inadmissibility for Case 1 when $a_1 = 0$ and $a_2 < 1 + k/2 - 1/k$ and, for the special case of Case 2 of known $\boldsymbol{\beta}$, when $a_1 = 0$ and $a_2 < (k + 1)/2 - 1/k$. Since the constant prior on $\mathbf{V}$ is $a_1 = a_2 = 0$, this clearly shows that the constant prior is badly inadmissible (i.e., is far from the boundary of admissibility). Kubokawa's results do not settle the question of admissibility of the hierarchical independence Jeffreys prior.

Either the constant prior or the hierarchical independence Jeffreys prior is easy to handle computationally so, if computational ease is the primary concern, our recommendation would be to use the hierarchical independence Jeffreys prior. As mentioned earlier, however, it is not immediately obvious how to generalize this prior to other hierarchical settings, although replacing $\mathbf{I}$ by the covariance matrix from the lower level is a good general solution when the lower level has an exchangeable structure.

The two proposed hierarchical reference priors, (a) $\pi(\mathbf{V}) = [|\mathbf{I} + \mathbf{V}| \prod_{i<j}(d_i - d_j)]^{-1}$ and (b) $\pi(\mathbf{V}) = [|\mathbf{V}|^{-(2k-1)/(2k)} \prod_{i<j}(d_i - d_j)]^{-1}$, are very appealing. They always result in proper posteriors if $m \geq 2$, a practically very useful and surprising fact when $m < k$ (explained in Section 2.2.4). They also both

yield admissible (or nearly admissible) estimators in Cases 2 and 3, and are computationally of similar complexity. Choice (a) is an actual hierarchical reference prior, in that it can be derived by a reference prior argument. In contrast, (b) was a rather ad hoc modification. Hence (a) should be the preferred choice for the actual model we consider. Again, however, it can be difficult in more general hierarchical models to know what to use in place of $\mathbf{I}$, and choice (b) does not require this additional input.

A very useful generalization (e.g., in common meta-analysis situations) would be to consider the setting

$$\mathbf{X}_i \sim N_k(\boldsymbol{\theta}_i, \boldsymbol{\Sigma}_i), \qquad \boldsymbol{\theta}_i \sim N_k(\mathbf{z}_i\boldsymbol{\beta}, \mathbf{V}),$$

independently for $i = 1, \ldots, m$, where the $\boldsymbol{\Sigma}_i$ are known positive-definite matrices, the $\mathbf{z}_i$ are given $k \times h$ covariate matrices and $\boldsymbol{\beta}$ is now $h \times 1$. Reasonable adaptations of the priors discussed above are:

1. Replace the covariance matrix $\mathbf{I}$ in the definitions of the priors for $\mathbf{V}$ by $\tilde{\boldsymbol{\Sigma}} = \frac{1}{m}\sum_{i=1}^m \boldsymbol{\Sigma}_i$. (Again, this is not necessary if one uses the prior $[|\mathbf{V}|^{-(2k-1)/(2k)} \prod_{i<j}(d_i - d_j)]^{-1}$.)
2. Replace the prior in (5) by $\pi(\boldsymbol{\beta}) = [1 + \boldsymbol{\beta}^t\mathbf{Z}^t\mathbf{Z}\boldsymbol{\beta}]^{-(h-1)/2}$, where $\mathbf{Z}$ is the matrix $(\mathbf{z}_1^t\mathbf{z}_2^t \cdots \mathbf{z}_m^t)^t$.

The results in the paper almost certainly go through for the generalization to known $\boldsymbol{\Sigma}_i$. We would also guess that the results are true for the generalization to covariates (the extension was true for the case $k = 1$, as shown in [7]), but the technical details in establishing this appear to be formidable. Finally, a number of the computational strategies mentioned in Section 1.5 are adaptable to these generalizations, but we do not have experience in utilization of such adaptations and so cannot comment on their efficiency.

## 2. Posterior propriety and impropriety.

2.1. *The marginal distribution.* Posterior propriety and admissibility properties are determined by study of the marginal density of $\mathbf{X}$, given by

$$m(\mathbf{x}) = \iiint f(\mathbf{x}|\boldsymbol{\theta})\pi(\boldsymbol{\theta}|\boldsymbol{\beta}, \mathbf{V})\pi(\boldsymbol{\beta})\pi(\mathbf{V})\, d\mathbf{V}\, d\boldsymbol{\beta}\, d\boldsymbol{\theta}, \tag{7}$$

where

$$\begin{aligned} f(\mathbf{x}|\boldsymbol{\theta}) &= \frac{1}{(2\pi)^{p/2}} \exp\left(-\frac{1}{2}(\mathbf{x} - \boldsymbol{\theta})^t(\mathbf{x} - \boldsymbol{\theta})\right) \\ &= \frac{1}{(2\pi)^{p/2}} \exp\left(-\frac{1}{2}\sum_{i=1}^m (\mathbf{x}_i - \boldsymbol{\theta}_i)^t(\mathbf{x}_i - \boldsymbol{\theta}_i)\right), \end{aligned} \tag{8}$$

$$
\begin{aligned}
\pi(\boldsymbol{\theta}|\boldsymbol{\beta},\mathbf{H},\mathbf{D}) &= \frac{1}{(2\pi)^{p/2}|\mathbf{V}|^{m/2}}\exp\left(-\frac{1}{2}\sum_{i=1}^{m}(\boldsymbol{\theta}_i-\boldsymbol{\beta})^t\mathbf{V}^{-1}(\boldsymbol{\theta}_i-\boldsymbol{\beta})\right)\\
&= \frac{1}{(2\pi)^{p/2}|\mathbf{D}|^{m/2}}\exp\left(-\frac{1}{2}\sum_{i=1}^{m}(\boldsymbol{\theta}_i-\boldsymbol{\beta})^t\mathbf{H}^t\mathbf{D}^{-1}\mathbf{H}(\boldsymbol{\theta}_i-\boldsymbol{\beta})\right).
\end{aligned} \tag{9}
$$

NOTATIONAL CONVENTION. It will be useful to write

$$m(\mathbf{x}) \approx g(\mathbf{x}) \tag{10}$$

if there exist $C_1 > 0$ and $C_2 > 0$ such that, $\forall\,\mathbf{x}$,

$$C_1 g(\mathbf{x}) \le m(\mathbf{x}) \le C_2 g(\mathbf{x}). \tag{11}$$

(This is related to the notion of "credence," as defined in [29].) Thus, under Condition 1 we can write

$$
\begin{aligned}
m(\mathbf{x}) \approx \iiiint & f(\mathbf{x}|\boldsymbol{\theta})\pi(\boldsymbol{\theta}|\boldsymbol{\beta},\mathbf{H},\mathbf{D})\pi(\boldsymbol{\beta})\\
&\times \frac{[\prod_{i<j}(d_i-d_j)]^l I_{[d_1>d_2>\cdots>d_k]}}{|\mathbf{I}+\mathbf{D}|^{(a_2-a_1)}|\mathbf{D}|^{a_1}}\,d\mathbf{D}\,d\mathbf{H}\,d\boldsymbol{\beta}\,d\boldsymbol{\theta}.
\end{aligned} \tag{12}
$$

Standard calculations yield the following expressions for $m(\mathbf{x})$ for the various cases of $\pi(\boldsymbol{\beta})$, where we define $\bar{\mathbf{x}} = \frac{1}{m}\sum_{i=1}^{m}\mathbf{x}_i$.

LEMMA 2.1. *For $\pi(\boldsymbol{\beta}) = 1$ (Case 1 scenario) and $m \ge 2$, the marginal density of $\mathbf{X}$ satisfies*

$$
\begin{aligned}
m(\mathbf{x}) \propto \iint & \frac{1}{|\mathbf{I}+\mathbf{D}|^{(m-1)/2}}\\
&\times \exp\left(-\frac{1}{2}\sum_{i=1}^{m}(\mathbf{x}_i-\bar{\mathbf{x}})^t\mathbf{H}^t(\mathbf{I}+\mathbf{D})^{-1}\mathbf{H}(\mathbf{x}_i-\bar{\mathbf{x}})\right)\pi(\mathbf{H},\mathbf{D})\,d\mathbf{D}\,d\mathbf{H}\\
\approx \iint & \frac{[\prod_{i<j}(d_i-d_j)]^l I_{[d_1>d_2>\cdots>d_k]}}{|\mathbf{I}+\mathbf{D}|^{[a_2-a_1+(m-1)/2]}|\mathbf{D}|^{a_1}}\\
&\times \exp\left(-\frac{1}{2}\sum_{i=1}^{m}(\mathbf{x}_i-\bar{\mathbf{x}})^t\mathbf{H}^t(\mathbf{I}+\mathbf{D})^{-1}\mathbf{H}(\mathbf{x}_i-\bar{\mathbf{x}})\right)d\mathbf{D}\,d\mathbf{H}.
\end{aligned} \tag{13}
$$

*When $m = 1$, the marginal density of $\mathbf{X}$ does not exist if $\pi(\mathbf{D})$ has infinite mass.*

Lemma 2.2. *If* $\pi(\boldsymbol{\beta})$ *is* $N_k(\mathbf{0},\mathbf{A})$ *(the Case* 2 *scenario, where we set* $\boldsymbol{\beta}^0=\mathbf{0}$ *for convenience), the marginal density of* $\mathbf{X}$ *is*

$$\begin{aligned}
m(\mathbf{x}) &\propto \iint \frac{1}{|\mathbf{I}+\mathbf{D}|^{(m-1)/2}|\mathbf{I}+\mathbf{D}+m\mathbf{H}\mathbf{A}\mathbf{H}^t|^{1/2}} \\
&\qquad \times \exp\left(-\frac{1}{2}\sum_{i=1}^m(\mathbf{x}_i-\bar{\mathbf{x}})^t\mathbf{H}^t(\mathbf{I}+\mathbf{D})^{-1}\mathbf{H}(\mathbf{x}_i-\bar{\mathbf{x}})\right) \\
&\qquad \times \exp\left(-\frac{1}{2}m(\mathbf{H}\bar{\mathbf{x}})^t(\mathbf{I}+\mathbf{D}+m\mathbf{H}\mathbf{A}\mathbf{H}^t)^{-1}(\mathbf{H}\bar{\mathbf{x}})\right)\pi(\mathbf{H},\mathbf{D})\,d\mathbf{D}\,d\mathbf{H} \\
&\approx \iint \frac{[\prod_{i<j}(d_i-d_j)]^l I_{[d_1>d_2>\cdots>d_k]}}{|\mathbf{I}+\mathbf{D}|^{[a_2-a_1+(m-1)/2]}|\mathbf{I}+\mathbf{D}+m\mathbf{H}\mathbf{A}\mathbf{H}^t|^{1/2}|\mathbf{D}|^{a_1}} \\
&\qquad \times \exp\left(-\frac{1}{2}\sum_{i=1}^m(\mathbf{x}_i-\bar{\mathbf{x}})^t\mathbf{H}^t(\mathbf{I}+\mathbf{D})^{-1}\mathbf{H}(\mathbf{x}_i-\bar{\mathbf{x}})\right) \\
&\qquad \times \exp\left(-\frac{1}{2}m(\mathbf{H}\bar{\mathbf{x}})^t(\mathbf{I}+\mathbf{D}+m\mathbf{H}\mathbf{A}\mathbf{H}^t)^{-1}(\mathbf{H}\bar{\mathbf{x}})\right)d\mathbf{D}\,d\mathbf{H}.
\end{aligned} \tag{14}$$

Lemma 2.3. *For* $\pi(\boldsymbol{\beta})=N_k(\mathbf{0},\lambda\mathbf{A})$ *(the Case* 3 *scenario, where we set* $\boldsymbol{\beta}^0=\mathbf{0}$ *for convenience), where* $\pi(\lambda)$ *satisfies Condition* 2*, the marginal density of* $\mathbf{X}$ *is*

$$\begin{aligned}
m(\mathbf{x}) &\propto \iiint \frac{1}{|\mathbf{I}+\mathbf{D}|^{(m-1)/2}|\mathbf{I}+\mathbf{D}+m\lambda\mathbf{H}\mathbf{A}\mathbf{H}^t|^{1/2}} \\
&\qquad \times \exp\left(-\frac{1}{2}\sum_{i=1}^m(\mathbf{x}_i-\bar{\mathbf{x}})^t\mathbf{H}^t(\mathbf{I}+\mathbf{D})^{-1}\mathbf{H}(\mathbf{x}_i-\bar{\mathbf{x}})\right) \\
&\qquad \times \exp\left(-\frac{1}{2}m(\mathbf{H}\bar{\mathbf{x}})^t(\mathbf{I}+\mathbf{D}+m\lambda\mathbf{H}\mathbf{A}\mathbf{H}^t)^{-1}(\mathbf{H}\bar{\mathbf{x}})\right) \\
&\qquad \times \pi(\mathbf{H},\mathbf{D})\pi(\lambda)\,d\lambda\,d\mathbf{D}\,d\mathbf{H} \\
&\approx \iiint \frac{[\prod_{i<j}(d_i-d_j)]^l I_{[d_1>d_2>\cdots>d_k]}}{|\mathbf{I}+\mathbf{D}|^{[a_2-a_1+(m-1)/2]}|\mathbf{I}+\mathbf{D}+m\lambda\mathbf{H}\mathbf{A}\mathbf{H}^t|^{1/2}|\mathbf{D}|^{a_1}} \\
&\qquad \times \exp\left(-\frac{1}{2}\sum_{i=1}^m(\mathbf{x}_i-\bar{\mathbf{x}})^t\mathbf{H}^t(\mathbf{I}+\mathbf{D})^{-1}\mathbf{H}(\mathbf{x}_i-\bar{\mathbf{x}})\right) \\
&\qquad \times \exp\left(-\frac{1}{2}m(\mathbf{H}\bar{\mathbf{x}})^t(\mathbf{I}+\mathbf{D}+m\lambda\mathbf{H}\mathbf{A}\mathbf{H}^t)^{-1}(\mathbf{H}\bar{\mathbf{x}})\right) \\
&\qquad \times \pi(\lambda)\,d\lambda\,d\mathbf{D}\,d\mathbf{H}.
\end{aligned} \tag{15}$$

2.2. *Impropriety of the posterior.* The next several theorems discuss the conditions under which the posterior distribution is proper. The following two lemmas are used.

LEMMA 2.4. *If a $k \times k$ matrix $\mathbf{H}$ is orthonormal, $\mathbf{h}_i(\mathbf{x})$, $i = 1, 2, \ldots, m$, are vector-valued functions, and $\mathbf{A}$ is positive semidefinite, then*

$$\begin{aligned} 0 &< \exp\Biggl(-\frac{1}{2}\sum_{i=1}^{m}\|\mathbf{h}_i(\mathbf{x})\|^2\Biggr) \\ &\le \exp\Biggl(-\frac{1}{2}\sum_{i=1}^{m}(\mathbf{h}_i(\mathbf{x}))^t\mathbf{H}^t(\mathbf{I}+\mathbf{D}+\mathbf{A})^{-1}\mathbf{H}\mathbf{h}_i(\mathbf{x})\Biggr) \le 1. \end{aligned} \tag{16}$$

PROOF. The upper bound is clear. On the other hand, $\mathbf{H}^t(\mathbf{I}+\mathbf{D}+\mathbf{A})^{-1}\mathbf{H} \le \mathbf{I}$ since $\mathbf{H}$ is orthonormal and $d_i \ge 0$, so that

$$\|(\mathbf{h}_i(\mathbf{x}))^t\mathbf{H}^t(\mathbf{I}+\mathbf{D}+\mathbf{A})^{-1}\mathbf{H}\mathbf{h}_i(\mathbf{x})\| \le \|\mathbf{h}_i(\mathbf{x})\|^2,$$

which yields the lower bound in (16). $\square$

LEMMA 2.5. *Let $\rho_1$ and $\rho_k$ be the maximum and the minimum eigenvalue of $\mathbf{A}$, respectively. Then*

$$\begin{aligned} |\mathbf{I}+\mathbf{D}| &\le |\mathbf{I}+\mathbf{D}+m\rho_k\mathbf{I}| \le |\mathbf{I}+\mathbf{D}+m\mathbf{H}\mathbf{A}\mathbf{H}^t| \\ &\le |\mathbf{I}+\mathbf{D}+m\rho_1\mathbf{I}| \le (1+m\rho_1)^k|\mathbf{I}+\mathbf{D}| \end{aligned} \tag{17}$$

*and*

$$\begin{aligned} \mathbf{x}^t(\mathbf{I}+\mathbf{D}+m\rho_1\mathbf{I})^{-1}\mathbf{x} &\le \mathbf{x}^t(\mathbf{I}+\mathbf{D}+m\mathbf{H}\mathbf{A}\mathbf{H}^t)^{-1}\mathbf{x} \\ &\le \mathbf{x}^t(\mathbf{I}+\mathbf{D}+m\rho_k\mathbf{I})^{-1}\mathbf{x}. \end{aligned} \tag{18}$$

PROOF. Using the notation $A \le B$ to denote that $B - A$ is nonnegative definite, we have

$$\rho_k\mathbf{I} \le \mathbf{H}\mathbf{A}\mathbf{H}^t \le \rho_1\mathbf{I}, \tag{19}$$

since $\mathbf{A}$ is nonnegative definite and $\mathbf{H}$ is orthonormal. Hence,

$$\mathbf{I}+\mathbf{D}+m\rho_k\mathbf{I} \le \mathbf{I}+\mathbf{D}+m\mathbf{H}\mathbf{A}\mathbf{H}^t \le \mathbf{I}+\mathbf{D}+m\rho_1\mathbf{I}, \tag{20}$$

from which (17) follows directly. From (20), clearly,

$$(\mathbf{I}+\mathbf{D}+m\rho_1\mathbf{I})^{-1} \le (\mathbf{I}+\mathbf{D}+m\mathbf{H}\mathbf{A}\mathbf{H}^t)^{-1} \le (\mathbf{I}+\mathbf{D}+m\rho_k\mathbf{I})^{-1}.$$

Equation (18) follows immediately, completing the proof. $\square$

Now we give the conditions under which the posterior distribution is proper for each of the three cases of $\pi(\boldsymbol{\beta})$.

2.2.1. *Case* 1 *scenario.* Since we are only considering improper $\pi(\mathbf{V})$, Lemma 2.1 shows that we need to consider only $m \geq 2$.

THEOREM 2.6. *If* $\pi(\boldsymbol{\beta}) = 1$, $m \geq 2$, $k \geq 2$, *and* $\pi(\mathbf{H}, \mathbf{D})$ *satisfies Condition* 1, *then the posterior distribution exists if and only if* $a_1 < 1$ *and* $a_2 > \frac{3-m}{2} + (k-1)l$.

PROOF. The posterior distribution is proper if and only if $0 < m(\mathbf{x}) < \infty$. The lower bound is clearly satisfied, so we only need to consider the upper bound. From (13) and Lemma 2.4, it is clear that, with $\mathbf{x}$ considered fixed, the posterior exists if and only if

$$m(\mathbf{x}) \approx \int \frac{[\prod_{i<j}(d_i - d_j)]^l I_{[d_1 > d_2 > \cdots > d_k]}}{|\mathbf{I} + \mathbf{D}|^{[a_2 - a_1 + (m-1)/2]} |\mathbf{D}|^{a_1}} \, d\mathbf{D} < \infty. \tag{21}$$

To determine necessary conditions for (21) to hold, first fix $d_1, d_2, \ldots, d_{k-1}$ and consider the integral over $d_k$ in (21), which is

$$C \int_0^{d_{k-1}} \frac{1}{d_k^{a_1}(1 + d_k)^{[a_2 - a_1 + (m-1)/2]}} \cdot \left[\prod_{i=1}^{k-1}(d_i - d_k)\right]^l dd_k.$$

Clearly,

$$\frac{1}{d_k^{a_1}(1 + d_k)^{[a_2 - a_1 + (m-1)/2]}} \cdot \left[\prod_{i=1}^{k-1}(d_i - d_k)\right]^l \sim \frac{C}{d_k^{a_1}} \qquad \text{as } d_k \to 0$$

and, when $a_1 \geq 1$,

$$\int_0^{d_{k-1}} \frac{1}{d_k^{a_1}} \, dd_k = \infty.$$

It follows that a necessary condition for (21) to hold is $a_1 < 1$.

Next, fix $d_2, d_3, \ldots, d_k$ and consider the integral over $d_1$ in (21),

$$C \int_{d_2}^{\infty} \frac{1}{d_1^{a_1}(1 + d_1)^{[a_2 - a_1 + (m-1)/2]}} \cdot \left[\prod_{i=2}^{k}(d_1 - d_i)\right]^l dd_1.$$

Counting the orders of $d_1$ for both the numerator and the denominator in the integral above, we see that this integral is infinite when $(k-1)l - (a_2 + (m-1)/2) \geq -1$. Thus another necessary condition for (21) to hold is

$$a_2 > \frac{3 - m}{2} + (k-1)l.$$

Next we show that the conditions given in the theorem are sufficient. Since $0 \le l \le 1$,

$$\int \frac{[\prod_{i<j}(d_i - d_j)]^l I_{[d_1>d_2>\cdots>d_k]}}{|\mathbf{I}+\mathbf{D}|^{[a_2-a_1+(m-1)/2]}|\mathbf{D}|^{a_1}}\,d\mathbf{D}$$

$$\le \int \frac{[\prod_{i<j}(d_i)]^l}{|\mathbf{I}+\mathbf{D}|^{[a_2-a_1+(m-1)/2]}|\mathbf{D}|^{a_1}}\,d\mathbf{D}$$

$$\le \prod_{i=1}^{k}\int_0^\infty \frac{d_i^{(k-i)l-a_1}}{(1+d_i)^{[a_2-a_1+(m-1)/2]}}\,dd_i.$$

Since $a_1 < 1$ and $(k-i)l \ge 0$, it is clear that each of these integrals is finite near 0. For $d_i$ near infinity, the corresponding integral is finite if $(3-m)/2 + (k-i)l < a_2$. This is true for all $i$ under the condition of the theorem, completing the proof. $\square$

2.2.2. *Case 2 scenario.*

THEOREM 2.7. *If $\boldsymbol{\beta} \sim N_k(\mathbf{0}, \mathbf{A})$, $k \ge 2$, and $\pi(\mathbf{H}, \mathbf{D})$ satisfies Condition 1, then the posterior distribution exists if and only if $a_1 < 1$ and $a_2 > 1 - \frac{m}{2} + (k-1)l$.*

PROOF. Clearly, we only need to find the necessary and sufficient condition for $m(\mathbf{x}) < \infty$. From (14) in Lemma 2.2 and Lemma 2.4, it is clear that

$$m(\mathbf{x}) \approx \iint \frac{[\prod_{i<j}(d_i - d_j)]^l I_{[d_1>d_2>\cdots>d_k]}}{|\mathbf{I}+\mathbf{D}|^{[a_2-a_1+(m-1)/2]}|\mathbf{I}+\mathbf{D}+m\mathbf{H}\mathbf{A}\mathbf{H}^t|^{1/2}|\mathbf{D}|^{a_1}}\,d\mathbf{D}\,d\mathbf{H}.$$

Again letting $\rho_1$ and $\rho_k$ denote the maximum and minimum eigenvalue of $\mathbf{A}$, it follows from (17) that $m(\mathbf{x}) < \infty$ if and only if

$$\int \frac{[\prod_{i<j}(d_i - d_j)]^l}{|\mathbf{I}+\mathbf{D}|^{(a_2-a_1+m/2)}|\mathbf{D}|^{a_1}} \cdot I_{[d_1>d_2>\cdots>d_k]} \cdot d\mathbf{D} < \infty.$$

The proof then proceeds in identical fashion to that of Theorem 2.6. $\square$

2.2.3. *Case 3 scenario.*

THEOREM 2.8. *Suppose that $\boldsymbol{\beta} \sim N_k(\mathbf{0}, \lambda\mathbf{A})$, $k \ge 2$, $\pi(\lambda)$ satisfies Condition 2 and $\pi(\mathbf{H}, \mathbf{D})$ satisfies Condition 1. The necessary and sufficient conditions for the posterior distribution to exist are $a_1 < 1$, $a_2 > 1 - \frac{m}{2} + (k-1)l$ and $b > 1 - \frac{k}{2}$.*

PROOF. As in the proof of Theorem 2.7, it is clear that

$$m(\mathbf{x}) \approx \iiint \frac{[\prod_{i<j}(d_i-d_j)]^l \cdot I_{[d_1>d_2>\cdots>d_k]}}{|\mathbf{I}+\mathbf{D}|^{(a_2-a_1+(m-1)/2)}|\mathbf{I}+\mathbf{D}+m\lambda\mathbf{H}\mathbf{A}\mathbf{H}^t|^{1/2}|\mathbf{D}|^{a_1}}\pi(\lambda)\,d\lambda\,d\mathbf{D}\,d\mathbf{H}.$$

Using (17), it is clear that the posterior density exists if and only if

$$I = \iint \frac{[\prod_{i<j}(d_i-d_j)]^l \cdot I_{[d_1>d_2>\cdots>d_k]}}{|\mathbf{I}+\mathbf{D}|^{(a_2-a_1+(m-1)/2)}|\mathbf{D}|^{a_1}|(1+c\lambda)\mathbf{I}+\mathbf{D}|^{1/2}}\pi(\lambda)\,d\lambda\,d\mathbf{D} < \infty. \tag{22}$$

Clearly,

$$\begin{aligned} I &\ge \int_0^1 \left\{\int \frac{[\prod_{i<j}(d_i-d_j)]^l I_{[d_1>d_2>\cdots>d_k]}}{|\mathbf{I}+\mathbf{D}|^{[a_2-a_1+(m-1)/2]}|\mathbf{D}|^{a_1}|(1+c\lambda)\mathbf{I}+\mathbf{D}|^{1/2}}\,d\mathbf{D}\right\}\pi(\lambda)\,d\lambda \\ &\ge \int_0^1 \left\{\int \frac{[\prod_{i<j}(d_i-d_j)]^l I_{[d_1>d_2>\cdots>d_k]}}{|\mathbf{I}+\mathbf{D}|^{[a_2-a_1+(m-1)/2]}|\mathbf{D}|^{a_1}|(1+c)\mathbf{I}+\mathbf{D}|^{1/2}}\,d\mathbf{D}\right\}\pi(\lambda)\,d\lambda \\ &\ge C\int \frac{[\prod_{i<j}(d_i-d_j)]^l I_{[d_1>d_2>\cdots>d_k]}}{|\mathbf{I}+\mathbf{D}|^{(a_2-a_1+m/2)}|\mathbf{D}|^{a_1}}\,d\mathbf{D}, \end{aligned}$$

the last inequality holding because of Condition 2(i). Proceeding as in Theorem 2.6, a necessary condition for $I$ to be finite is

$$a_1 < 1 \quad \text{and} \quad a_2 > 1 - \frac{m}{2} + (k-1)l.$$

On the other hand, by (ii) of Condition 2,

$$\begin{aligned} I &\ge \int_1^\infty \left\{\int \frac{[\prod_{i<j}(d_i-d_j)]^l I_{[d_1>d_2>\cdots>d_k]}}{|\mathbf{I}+\mathbf{D}|^{[a_2-a_1+(m-1)/2]}|\mathbf{D}|^{a_1}|(1+c\lambda)\mathbf{I}+\mathbf{D}|^{1/2}}\,d\mathbf{D}\right\}\pi(\lambda)\,d\lambda \\ &\ge C\int_1^\infty \frac{1}{(1+c\lambda)^{k/2}}\cdot\frac{1}{\lambda^b}\,d\lambda \cdot \int \frac{[\prod_{i<j}(d_i-d_j)]^l I_{[d_1>d_2>\cdots>d_k]}}{|\mathbf{I}+\mathbf{D}|^{(a_2-a_1+m/2)}|\mathbf{D}|^{a_1}}\,d\mathbf{D}. \end{aligned}$$

This integral is infinite when $b \le 1-k/2$. So another necessary condition for (22) to hold is $b > 1-k/2$.

Next let us prove that the conditions are sufficient. Using

$$\prod_{j=1}^k \frac{1}{(1+C\lambda+d_j)^{1/2}} \le \frac{1}{(1+C\lambda+d_1)^{1/2}(1+C\lambda)^{(k-1)/2}},$$

we have

$$\begin{aligned} I &\le \iint \frac{[\prod_{i<j}(d_i-d_j)]^l I_{[d_1>d_2>\cdots>d_k]}}{|\mathbf{I}+\mathbf{D}|^{[a_2-a_1+(m-1)/2]}(1+C\lambda+d_1)^{1/2}(1+C\lambda)^{(k-1)/2}|\mathbf{D}|^{a_1}}\pi(\lambda)\,d\lambda\,d\mathbf{D} \\ &\le \iint \frac{(\prod_{i<j} d_i)^l}{|\mathbf{I}+\mathbf{D}|^{[a_2-a_1+(m-1)/2]}(1+C\lambda+d_1)^{1/2}(1+C\lambda)^{(k-1)/2}|\mathbf{D}|^{a_1}}\pi(\lambda)\,d\lambda\,d\mathbf{D}. \end{aligned}$$

As in the proof of Theorem 2.7, the integrals over $d_2 \cdots d_p$ are finite under the stated conditions, so that

$$I \leq C \iint \frac{d_1^{[(k-1)l-a_1]}}{(1+d_1)^{[a_2-a_1+(m-1)/2]}(1+C\lambda+d_1)^{1/2}(1+C\lambda)^{(k-1)/2}} \pi(\lambda)\, d\lambda\, dd_1$$
$$\leq C \iint (1+d_1)^{-[a_2+(m-1)/2-(k-1)l]}$$
$$\times (1+C\lambda+d_1)^{-1/2}(1+C\lambda)^{-(k-1)/2}\pi(\lambda)\, d\lambda\, dd_1.$$

Break this integral up into four integrals over $(0,c)\times(0,c)$, $(0,c)\times(c,\infty)$, $(c,\infty)\times(0,c)$ and $(c,\infty)\times(c,\infty)$. Bounding the first three integrals is easy, using Condition 2. The last integral is bounded as in the proof of Lemma 1 of [7].

□

2.2.4. *Summary of posterior propriety and impropriety.* The cases of most interest are $l=0$ and $l=1$. The following corollaries of Theorems 2.6, 2.7 and 2.8 deal with these cases.

COROLLARY 2.9. *Suppose $l=0$ and $k\geq 2$.*

(a) *In the Case 1 scenario ($\pi(\boldsymbol{\beta})=1$), when $m\geq 2$, the posterior distribution exists if and only if $a_1<1$ and $a_2>\frac{3-m}{2}$.*

(b) *In the Case 2 scenario ($\boldsymbol{\beta}\sim N_k(\mathbf{0},\mathbf{A})$), the posterior distribution exists if and only if $a_1<1$ and $a_2>1-\frac{m}{2}$.*

(c) *In the Case 3 scenario ($\boldsymbol{\beta}\sim N_k(\mathbf{0},\lambda\mathbf{A})$, $\lambda\sim\pi(\lambda)$), the posterior distribution exists if and only if $a_1<1$, $a_2>1-\frac{m}{2}$ and $b>1-\frac{k}{2}$.*

COROLLARY 2.10. *Suppose $l=1$ and $k\geq 2$.*

(a) *In the Case 1 scenario ($\pi(\boldsymbol{\beta})=1$), when $m\geq 2$, the posterior distribution exists if and only if $a_1<1$ and $a_2>k-\frac{m-1}{2}$.*

(b) *In the Case 2 scenario ($\boldsymbol{\beta}\sim N_k(\mathbf{0},\mathbf{A})$), the posterior distribution exists if and only if $a_1<1$ and $a_2>k-\frac{m}{2}$.*

(c) *In the Case 3 scenario ($\boldsymbol{\beta}\sim N_k(\mathbf{0},\lambda\mathbf{A})$, $\lambda\sim\pi(\lambda)$), the posterior distribution exists when $a_1<1$, $a_2>k-\frac{m}{2}$ and $b>1-\frac{k}{2}$.*

It follows that the most commonly used objective priors for covariance matrices cannot be used in the hierarchical setting. The nonhierarchical independence Jeffreys prior $[l=1, a_1=a_2=(k+1)/2]$ and the nonhierarchical reference prior $(l=0, a_1=a_2=1)$ yield improper posteriors. The constant prior $(l=1, a_1=a_2=0)$ yields a proper posterior only when $2k<m-1$ for Case 1, and when $2k<m$ for Case 2 and Case 3. This implies that the number of blocks $m$ has to be at least $2k+2$ for Case 1 and $2k+1$ for Case 2 and Case 3.

In contrast, the hierarchical independence Jeffreys prior $[l=1, a_1=0, a_2=(k+1)/2]$ yields a proper posterior when $m>k$ for Case 1 and $m>k-1$ for Cases 2 and 3, considerably weaker conditions. Furthermore, the hierarchical reference prior (a) $(l=0,\ a_1=0,\ a_2=1)$ and the hierarchical reference prior (b) $[l=0,\ a_1=a_2=(2k-1)/(2k)]$ always yield a proper posterior, except when $m=1$ in Case 1.

It is quite surprising that posterior propriety for the hierarchical reference priors does not require $m$ to grow with $k$ (as is necessary for the hierarchical independence Jeffreys prior). One needs on the order of $m=k$ blocks in order for the hyper-covariance matrix $\mathbf{V}$ to be identifiable, which is usually viewed as being equivalent to posterior propriety. Such equivalence is clearly not the case here; in the simplest Case 1 scenario, for instance, only $m=2$ blocks are needed for posterior propriety of the reference priors, regardless of the value of $k$.

To understand why this is so, consider the transformed version of the problem in Section 1.4. Note that the domain of $\mathbf{H}$ is a compact set and the reference prior assigns a proper uniform distribution to this set, so the only parameters that intuitively need data to have proper posteriors are $\boldsymbol{\beta}$ and $\mathbf{D}$. These vectors consist of $2k$ unknowns, which intuitively can be handled by the $2k$ coordinate observations corresponding to $m=2$. This general posterior propriety is a very attractive property of the hierarchical reference priors in that it is often difficult in complicated hierarchical models to ensure that conditions such as $m>k$ are satisfied at all levels and components of the hierarchy.

## 3. Admissibility and inadmissibility.

3.1. *Introduction.* In this section we give conditions under which the hierarchical Bayes estimate $\boldsymbol{\delta}^{\pi}(\mathbf{x})$ (the posterior mean) of $\boldsymbol{\theta}$ is admissible and inadmissible for quadratic loss (2). We restrict consideration to the priors for which $l=0$, since these are the priors we will recommend and analysis for $l>0$ requires different techniques.

Our study utilizes the following powerful results from [9]. Define

$$\overline{m}(r)=\int m(\mathbf{x})\,d\phi(\mathbf{x}), \tag{23}$$

$$\underline{m}(r)=\int \frac{1}{m(\mathbf{x})}\,d\phi(\mathbf{x}), \tag{24}$$

where $\phi(\cdot)$ is the uniform probability measure on the surface of the sphere of radius $r=\|\mathbf{x}\|$.

RESULT 3.1. If $\boldsymbol{\delta}^\pi(\mathbf{x}) - \mathbf{x}$ is uniformly bounded and

$$\int_c^\infty [r^{mk-1}\overline{m}(r)]^{-1}\,dr = \infty \tag{25}$$

for some $c > 0$, then $\boldsymbol{\delta}^\pi(\mathbf{x})$ is admissible.

RESULT 3.2. If

$$\int_c^\infty r^{1-mk}\underline{m}(r)\,dr < \infty \tag{26}$$

for some $c > 0$, then $\boldsymbol{\delta}^\pi(\mathbf{x})$ is inadmissible.

3.2. *Preliminary lemmas.* The following lemmas are needed.

LEMMA 3.3. (a) *If $a < 1$, $r + a > 1$ and $c_1$ and $c_2$ are positive constants, then*

$$f(v) \equiv \int_0^\infty \frac{1}{(c_1+d)^r d^a} \exp\left(-\frac{v}{2(c_2+d)}\right) dd \approx C_1 \min\{C_2, v^{1-r-a}\}, \tag{27}$$

*for some positive constants $C_1$ and $C_2$.*

(b) *If $a > -1$, $\mu > 0$ and $v > 0$, then*

$$g(\mu, v) \equiv \int_0^\mu t^a e^{-vt}\,dt \le C \min\{v^{-(a+1)}, \mu^{(a+1)}\} \tag{28}$$

*for some positive constant $C$.*

For the proof see the Appendix.

LEMMA 3.4. *Assuming the integrals exist,*

$$\iint g(\mathbf{H}^t\mathbf{D}\mathbf{H}) I_{[d_1>d_2>\cdots>d_k]}\,d\mathbf{D}\,d\mathbf{H} = \frac{1}{k!}\iint g(\mathbf{H}^t\mathbf{D}\mathbf{H})\,d\mathbf{D}\,d\mathbf{H}. \tag{29}$$

PROOF. Suppose that $d_1 > d_2 > \cdots > d_k > 0$ are eigenvalues of $\mathbf{V}$ and $(d_1^*, d_2^*, \ldots, d_k^*)$ is a different ordering of $(d_1, d_2, \ldots, d_k)$. Let $\mathbf{D}^* = \operatorname{diag}(d_1^*, d_2^*, \ldots, d_k^*)$. Since there exists an orthonormal matrix $\mathbf{H}^*$ such that $\mathbf{D} = \mathbf{H}^{*t}\mathbf{D}^*\mathbf{H}^*$, it follows that

$$\begin{aligned}
&\iint g(\mathbf{H}^t\mathbf{D}\mathbf{H}) I_{[d_1>d_2>\cdots>d_k]}\,d\mathbf{D}\,d\mathbf{H} \\
&\quad = \iint g((\mathbf{H}^*\mathbf{H})^t\mathbf{D}^*(\mathbf{H}^*\mathbf{H})) I_{[d_1>d_2>\cdots>d_k]}\,d\mathbf{D}\,d\mathbf{H} \\
&\quad = \iint g((\mathbf{H}^*\mathbf{H})^t\mathbf{D}^*(\mathbf{H}^*\mathbf{H})) I^*\,d\mathbf{D}^*\,d\mathbf{H},
\end{aligned} \tag{30}$$

the last step following from the change of variables from $\mathbf{D} \to \mathbf{D}^*$ (which has Jacobian 1), where $I^*$ corresponds to the new ordering. Next, note that, since $d\mathbf{H}$ represents the invariant Haar density,

$$\int\int g((\mathbf{H}^*\mathbf{H})^t\mathbf{D}^*(\mathbf{H}^*\mathbf{H}))I^*\,d\mathbf{D}^*\,d\mathbf{H} = \int\int g(\mathbf{H}^t\mathbf{D}^*\mathbf{H})I^*\,d\mathbf{D}^*\,d\mathbf{H}.$$

Hence $\int g(\mathbf{H}^t\mathbf{D}^*\mathbf{H})\tilde{I}\,d\mathbf{D}$ is the same for any ordering $\tilde{I}$ of the eigenvalues, and the result follows since there are $k!$ orderings. $\square$

NOTATIONAL CONVENTION. We need to generalize the notation in (10). Indeed, let

$$m(\mathbf{x}) \approx g(\mathbf{c},\mathbf{x}) \quad \text{stand for} \quad g(\mathbf{c},\mathbf{x}) \le m(\mathbf{x}) \le g(\mathbf{c}',\mathbf{x}) \tag{31}$$

for some (possibly vectors) $\mathbf{c}$ and $\mathbf{c}'$. For instance, in (33) below, $\mathbf{c} = (C_1, C_2, C_3)$. The earlier notation was the special case where $g(c,\mathbf{x}) = cg(\mathbf{x})$.

We conclude this section with presentation of needed upper and lower bounds (using the $\approx$ notation) for the marginal densities in Cases 1, 2 and 3.

LEMMA 3.5. *In the Case* 1 *scenario and with* $l = 0$,

$$\begin{aligned} m(\mathbf{x}) \approx C \int\int & \frac{1}{|\mathbf{D}|^{a_1}|\mathbf{I}+\mathbf{D}|^{(m-1)/2+a_2-a_1}} \\ & \times \exp\left(-\frac{1}{2}\sum_{i=1}^m (\mathbf{x}_i - \bar{\mathbf{x}})^t\mathbf{H}^t(\mathbf{I}+\mathbf{D})^{-1}\mathbf{H}(\mathbf{x}_i - \bar{\mathbf{x}})\right)\cdot d\mathbf{D}\,d\mathbf{H}. \end{aligned} \tag{32}$$

PROOF. This follows directly from (13) in Lemma 2.1 and Lemma 3.4. $\square$

LEMMA 3.6. *In the Case* 2 *scenario and with* $l = 0$,

$$\begin{aligned} m(\mathbf{x}) \approx C_1 \int\int & \frac{1}{|\mathbf{D}|^{a_1}|C_2\mathbf{I}+\mathbf{D}|^{m/2+a_2-a_1}} \\ & \times \exp\left(-\frac{1}{2}\sum_{i=1}^m \mathbf{x}_i^t\mathbf{H}^t(C_3\mathbf{I}+\mathbf{D})^{-1}\mathbf{H}\mathbf{x}_i\right)\cdot d\mathbf{D}\,d\mathbf{H}. \end{aligned} \tag{33}$$

PROOF. From (14) in Lemma 2.2 and Lemma 3.4,

$$\begin{aligned} m(\mathbf{x}) \approx \int\int & \frac{1}{|\mathbf{D}|^{a_1}|\mathbf{I}+\mathbf{D}|^{(m-1)/2+a_2-a_1}|\mathbf{I}+\mathbf{D}+m\mathbf{H}\mathbf{A}\mathbf{H}^t|^{1/2}} \\ & \times \exp\left(-\frac{1}{2}\sum_{i=1}^m (\mathbf{x}_i - \bar{\mathbf{x}})^t\mathbf{H}^t(\mathbf{I}+\mathbf{D})^{-1}\mathbf{H}(\mathbf{x}_i - \bar{\mathbf{x}})\right) \\ & \times \exp\left(-\frac{1}{2}m(\mathbf{H}\bar{\mathbf{x}})^t(\mathbf{I}+\mathbf{D}+m\mathbf{H}\mathbf{A}\mathbf{H}^t)^{-1}(\mathbf{H}\bar{\mathbf{x}})\right)\cdot d\mathbf{D}\,d\mathbf{H}. \end{aligned} \tag{34}$$

Applying Lemma 2.5 to (34),

$$
\begin{aligned}
m(\mathbf{x}) &\le C \iint \frac{1}{|\mathbf{D}|^{a_1}|\mathbf{I}+\mathbf{D}|^{(m-1)/2+a_2-a_1}|\mathbf{I}+\mathbf{D}+m\rho_k\mathbf{I}|^{1/2}} \\
&\qquad \times \exp\Biggl(-\frac{1}{2}\sum_{i=1}^{m}(\mathbf{x}_i-\bar{\mathbf{x}})^t\mathbf{H}^t(\mathbf{I}+\mathbf{D})^{-1}\mathbf{H}(\mathbf{x}_i-\bar{\mathbf{x}})\Biggr) \\
&\qquad \times \exp\biggl(-\frac{1}{2}m(\mathbf{H}\bar{\mathbf{x}})^t(\mathbf{I}+\mathbf{D}+m\rho_1\mathbf{I})^{-1}(\mathbf{H}\bar{\mathbf{x}})\biggr)\cdot d\mathbf{D}\,d\mathbf{H} \\
&\le C \iint \frac{1}{|\mathbf{D}|^{a_1}|\mathbf{I}+\mathbf{D}|^{(m-1)/2+a_2-a_1}|\mathbf{I}+\mathbf{D}|^{1/2}} \\
&\qquad \times \exp\Biggl(-\frac{1}{2}\sum_{i=1}^{m}(\mathbf{x}_i-\bar{\mathbf{x}})^t\mathbf{H}^t(\mathbf{I}+\mathbf{D}+m\rho_1\mathbf{I})^{-1}\mathbf{H}(\mathbf{x}_i-\bar{\mathbf{x}})\Biggr) \\
&\qquad \times \exp\biggl(-\frac{1}{2}m(\mathbf{H}\bar{\mathbf{x}})^t(\mathbf{I}+\mathbf{D}+m\rho_1\mathbf{I})^{-1}(\mathbf{H}\bar{\mathbf{x}})\biggr)\cdot d\mathbf{D}\,d\mathbf{H} \\
&\le C \iint \frac{1}{|\mathbf{D}|^{a_1}|\mathbf{I}+\mathbf{D}|^{m/2+a_2-a_1}} \\
&\qquad \times \exp\Biggl(-\frac{1}{2}\sum_{i=1}^{m}\mathbf{x}_i^t\mathbf{H}^t(\mathbf{I}+\mathbf{D}+m\rho_1\mathbf{I})^{-1}\mathbf{H}\mathbf{x}_i\Biggr)\cdot d\mathbf{D}\,d\mathbf{H}.
\end{aligned}
$$

Similarly,

$$
\begin{aligned}
m(\mathbf{x}) &\ge C \iint \frac{1}{|\mathbf{D}|^{a_1}|\mathbf{I}+\mathbf{D}+m\rho_1\mathbf{I}|^{m/2+a_2-a_1}} \\
&\qquad \times \exp\Biggl(-\frac{1}{2}\sum_{i=1}^{m}\mathbf{x}_i^t\mathbf{H}^t(\mathbf{I}+\mathbf{D})^{-1}\mathbf{H}\mathbf{x}_i\Biggr)\cdot d\mathbf{D}\,d\mathbf{H}.
\end{aligned}
$$

This completes the proof. □

LEMMA 3.7. *In the Case 3 scenario and with $l=0$,*

$$
\begin{aligned}
m(\mathbf{x}) &\approx C_1 \iiint \frac{1}{|\mathbf{D}|^{a_1}|\mathbf{I}+\mathbf{D}|^{(m-1)/2+a_2-a_1}|(1+C_2\lambda)\mathbf{I}+\mathbf{D}|^{1/2}} \\
&\qquad \times \exp\Biggl(-\frac{1}{2}\sum_{i=1}^{m}(\mathbf{x}_i-\bar{\mathbf{x}})^t\mathbf{H}^t(\mathbf{I}+\mathbf{D})^{-1}\mathbf{H}(\mathbf{x}_i-\bar{\mathbf{x}})\Biggr) \\
&\qquad \times \exp\biggl(-\frac{1}{2}m(\mathbf{H}\bar{\mathbf{x}})^t[(1+C_3\lambda)\mathbf{I}+\mathbf{D}]^{-1}(\mathbf{H}\bar{\mathbf{x}})\biggr) \\
&\qquad \times \pi(\lambda)\,d\lambda\,d\mathbf{D}\,d\mathbf{H}.
\end{aligned}
\tag{35}
$$

PROOF. From (15) in Lemma 2.3 and Lemma 3.4,

$$
m(\mathbf{x}) \approx \iiint \frac{1}{|\mathbf{D}|^{a_1}|\mathbf{I}+\mathbf{D}|^{(m-1)/2+a_2-a_1}|\mathbf{I}+\mathbf{D}+m\lambda\mathbf{H}\mathbf{A}\mathbf{H}^t|^{1/2}}
$$

$$\times \exp\left(-\frac{1}{2}\sum_{i=1}^{m}(\mathbf{x}_i-\bar{\mathbf{x}})^t\mathbf{H}^t(\mathbf{I}+\mathbf{D})^{-1}\mathbf{H}(\mathbf{x}_i-\bar{\mathbf{x}})\right)$$

$$\times \exp\left(-\frac{1}{2}m(\mathbf{H}\bar{\mathbf{x}})^t(\mathbf{I}+\mathbf{D}+m\lambda\mathbf{H}\mathbf{A}\mathbf{H}^t)^{-1}(\mathbf{H}\bar{\mathbf{x}})\right)\cdot\pi(\lambda)\,d\lambda\,d\mathbf{D}\,d\mathbf{H}.$$

The proof is then exactly like that of Lemma 3.6. □

3.3. *Uniformly bounded property.* Let $\boldsymbol{\delta}^\pi(\mathbf{x})$ be the posterior mean of $\boldsymbol{\theta}$ with respect to the posterior distribution. To prove admissibility by Brown's results, we first need to show that $\boldsymbol{\delta}^\pi(\mathbf{x})-\mathbf{x}$ is uniformly bounded. Let $\boldsymbol{\delta}^\pi(\mathbf{x})_{p\times 1}=(\boldsymbol{\delta}_1^\pi(\mathbf{x}),\boldsymbol{\delta}_2^\pi(\mathbf{x}),\dots,\boldsymbol{\delta}_m^\pi(\mathbf{x}))^t$, so that $\boldsymbol{\delta}_i^\pi(\mathbf{x})$ is the subvector of $\boldsymbol{\delta}^\pi(\mathbf{x})$ corresponding to $\boldsymbol{\theta}_i$. By symmetry, it is clearly sufficient to show that $\boldsymbol{\delta}_1^\pi(\mathbf{x})-\mathbf{x}_1$ is uniformly bounded.

LEMMA 3.8. *Suppose that $\mathbf{z}_1,\mathbf{z}_2,\dots,\mathbf{z}_m$ are $k\times 1$ vectors and $\tilde{\mathbf{y}}$ is the $k\times 1$ vector $(y_1,y_2,\dots,y_k)^t$, with $y_n=(\sum_{i=1}^m z_{in}^2)^{1/2}$, where $z_{ij}$ is the jth element of $\mathbf{z}_i$. Define*

$$g(\mathbf{c},\mathbf{D},\mathbf{z}_1,\mathbf{z}_2,\dots,\mathbf{z}_m)=\frac{C_1}{|\mathbf{D}|^u|C_2\mathbf{I}+\mathbf{D}|^v}\exp\left(-\frac{1}{2}\sum_{i=1}^m\mathbf{z}_i^t(C_3\mathbf{I}+\mathbf{D})^{-1}\mathbf{z}_i\right),$$

*where the $C_i$ are positive constants. If $u+v>1$ and $u<1$, then*

$$\frac{\int\|(\mathbf{I}+\mathbf{D})^{-1}\mathbf{y}\|g(\mathbf{c},\mathbf{D},\mathbf{z}_1,\mathbf{z}_2,\dots,\mathbf{z}_m)\,d\mathbf{D}}{\int g(\mathbf{c}',\mathbf{D},\mathbf{z}_1,\mathbf{z}_2,\dots,\mathbf{z}_m)\,d\mathbf{D}}\tag{36}$$

*is uniformly bounded over $\mathbf{z}_1,\mathbf{z}_2,\dots,\mathbf{z}_m$.*

PROOF. In (36),

$$\begin{aligned}|\text{Numerator}|&\le\int\|(\mathbf{I}+\mathbf{D})^{-1}\mathbf{y}\|\frac{C_1}{|\mathbf{D}|^u|C_2\mathbf{I}+\mathbf{D}|^v}\\&\qquad\times\exp\left(-\frac{1}{2}\sum_{i=1}^m\mathbf{z}_i^t(C_3\mathbf{I}+\mathbf{D})^{-1}\mathbf{z}_i\right)d\mathbf{D}\\&\le C_1\sum_{n=1}^k\int\frac{|y_n|}{1+d_n}\left(\prod_{j=1}^k\frac{1}{d_j^u(C_2+d_j)^v}\right)\\&\qquad\times\exp\left(-\frac{1}{2(C_3+d_j)}\sum_{i=1}^m z_{ij}^2\right)dd_j.\end{aligned}$$

For each $n$ we will bound the $k$-dimensional integral. If $j\ne n$, by Lemma 3.3 with $a=u$ and $r=v$ it is clear that

$$\int_0^\infty\frac{1}{d_j^u(C_2+d_j)^v}\exp\left(-\frac{1}{2(C_3+d_j)}\sum_{i=1}^m z_{ij}^2\right)dd_j$$

$$\leq C_1 \min\left\{C_2, \left(\sum_{i=1}^m z_{ij}^2\right)^{1-u-v}\right\}.$$

If $j = n$, applying Lemma 3.3 with $a = u$ and $r = v + 1$ [again using the fact that $(1 + d_n)/(C_2 + d_n)$ is uniformly bounded] yields

$$C_1 \int_0^\infty \frac{1}{1+d_j} \cdot \frac{1}{d_j^u (C_2 + d_j)^v} \exp\left(-\frac{1}{2(C_3 + d_j)} \sum_{i=1}^m z_{ij}^2\right) dd_j$$
$$\leq C_1' \min\left\{C_2', \left(\sum_{i=1}^m z_{ij}^2\right)^{-u-v}\right\}.$$

Therefore,

$$|\text{Numerator}|$$
$$\leq \sum_{n=1}^k \left[|y_n| C_1' \min\left\{C_2', \left(\sum_{i=1}^m z_{in}^2\right)^{-u-v}\right\} \prod_{j \neq n} C_1 \min\left\{C_2, \left(\sum_{i=1}^m z_{ij}^2\right)^{1-u-v}\right\}\right].$$

In (36),

$$\text{Denominator} = C_1' \prod_{j=1}^k \int_0^\infty \frac{1}{d_j^u (C_2' + d_j)^v} \exp\left(-\frac{1}{2(C_3' + d_j)} \sum_{i=1}^m z_{ij}^2\right) dd_j.$$

Applying the lower bound of Lemma 3.3, with $a = u$ and $r = v$, yields

$$C_1' \int_0^\infty \frac{1}{d_j^u (C_2' + d_j)^v} \exp\left(-\frac{1}{2(C_3' + d_j)} \sum_{i=1}^m z_{ij}^2\right) dd_j$$
$$\geq C_1^* \min\left\{C_2^*, \left(\sum_{i=1}^m z_{ij}^2\right)^{1-u-v}\right\}.$$

Thus

$$\text{Denominator} \geq \prod_{i=1}^k \left[C_1^* \min\left\{C_2^*, \left(\sum_{i=1}^m z_{ij}^2\right)^{1-u-v}\right\}\right].$$

Combining the numerator and the denominator, we have

$$\left|\frac{\int \mathbf{H}^t(\mathbf{I}+\mathbf{D})^{-1}\mathbf{y} g(\mathbf{c}, \mathbf{D}, \mathbf{z}_1, \mathbf{z}_2, \ldots, \mathbf{z}_m)\, d\mathbf{D}}{\int g(\mathbf{c}', \mathbf{D}, \mathbf{z}_1, \mathbf{z}_2, \ldots, \mathbf{z}_m)\, d\mathbf{D}}\right|$$
$$\leq \left(\sum_{n=1}^k \left[|y_n| C_1' \min\left\{C_2', \left(\sum_{i=1}^m z_{in}^2\right)^{-u-v}\right\}\right.\right.$$
$$\left.\left. \times \prod_{j \neq n} C_1 \min\left\{C_2, \left(\sum_{i=1}^m z_{ij}^2\right)^{1-u-v}\right\}\right]\right)$$

$$\times \left(\prod_{j=1}^{k}\left[C_1^* \min\left\{C_2^*, \left(\sum_{i=1}^{m} z_{ij}^2\right)^{1-u-v}\right\}\right]\right)^{-1}$$

$$= \sum_{n=1}^{k}\left(\frac{|y_n| C_1' \min\{C_2', (\sum_{i=1}^m z_{in}^2)^{-u-v}\}}{C_1^* \min\{C_2^*, (\sum_{i=1}^m z_{in}^2)^{1-u-v}\}} \cdot \prod_{j\neq n} \frac{C_1 \min\{C_2, (\sum_{i=1}^m z_{ij}^2)^{1-u-v}\}}{C_1^* \min\{C_2^*, (\sum_{i=1}^m z_{ij}^2)^{1-u-v}\}}\right).$$

Clearly

$$\prod_{j\neq n} \frac{C_1 \min\{C_2, (\sum_{i=1}^m z_{ij}^2)^{1-u-v}\}}{C_1^* \min\{C_2^*, (\sum_{i=1}^m z_{ij}^2)^{1-u-v}\}} \le C.$$

Using the condition $u+v>1$, we have that for large $\sum_{i=1}^m z_{in}^2 = y_n^2$

$$\frac{|y_n| C_1' \min\{C_2', (\sum_{i=1}^m z_{in}^2)^{-u-v}\}}{C_1^* \min\{C_2^*, (\sum_{i=1}^m z_{in}^2)^{1-u-v}\}}$$

behaves as $|y_n|/y_n^2$, while for small $y_n^2$ it behaves as $C_3|y_n| < C_4$, so that the ratio is clearly uniformly bounded. Thus

$$\frac{\int \|(\mathbf{I}+\mathbf{D})^{-1}\mathbf{y}\| g(\mathbf{c}, \mathbf{D}, \mathbf{z}_1, \mathbf{z}_2, \ldots, \mathbf{z}_m)\, d\mathbf{D}}{\int g(\mathbf{c}', \mathbf{D}, \mathbf{z}_1, \mathbf{z}_2, \ldots, \mathbf{z}_m)\, d\mathbf{D}} \le C,$$

completing the proof. □

THEOREM 3.9. *Assume that* $\pi(\boldsymbol{\beta}) = 1$, $m \ge 2$, $k \ge 2$, *and* $\pi(\mathbf{H}, \mathbf{D})$ *satisfies Condition 1. Also suppose that we choose* $l = 0$. *If* $a_1 < 1$ *and* $a_2 > \frac{3-m}{2}$, *then* $\boldsymbol{\delta}^\pi(\mathbf{x}) - \mathbf{x}$ *is uniformly bounded.*

PROOF. We only need to show that $\boldsymbol{\delta}_1^\pi(\mathbf{x}) - \mathbf{x}_1$ is uniformly bounded. It is well known that

$$\boldsymbol{\delta}_1^\pi(\mathbf{x}) - \mathbf{x}_1 = (\nabla m(\mathbf{x}))_1 / m(\mathbf{x}), \tag{37}$$

where $\nabla$ denotes the gradient. Exactly as in the proof of Lemma 3.5, it can be shown that

$$\begin{aligned}\|(\nabla m(\mathbf{x}))_1\| = \Bigg\| -\iint & \frac{(\mathbf{I}+\mathbf{H}^t\mathbf{D}\mathbf{H})^{-1}(\mathbf{x}_1 - \bar{\mathbf{x}})}{|\mathbf{I}+\mathbf{D}|^{(m-1)/2}} \\ & \times \exp\left(-\frac{1}{2}\sum_{i=1}^m (\mathbf{x}_i - \bar{\mathbf{x}})^t \mathbf{H}^t (\mathbf{I}+\mathbf{D})^{-1} \mathbf{H}(\mathbf{x}_i - \bar{\mathbf{x}})\right) \\ & \times \pi(\mathbf{H}, \mathbf{D})\, d\mathbf{H}\, d\mathbf{D} \Bigg\| \\ \le \iint & \frac{\|(\mathbf{I}+\mathbf{D})^{-1}\mathbf{H}(\mathbf{x}_1 - \bar{\mathbf{x}})\|}{|\mathbf{D}|^{a_1} |\mathbf{I}+\mathbf{D}|^{(m-1)/2 + a_2 - a_1}}\end{aligned}$$

$$\times \exp\left(-\frac{1}{2}\sum_{i=1}^{m}(\mathbf{x}_i - \bar{\mathbf{x}})^t \mathbf{H}^t (\mathbf{I}+\mathbf{D})^{-1}\mathbf{H}(\mathbf{x}_i - \bar{\mathbf{x}})\right) \cdot d\mathbf{D}\, d\mathbf{H}.$$

Hence, defining $\mathbf{z}_i = \mathbf{H}(\mathbf{x}_i - \bar{\mathbf{x}})$, and using the lower bound in Lemma 3.5 for the denominator in (37), one obtains (for appropriate constants $\mathbf{c}$, $\mathbf{c}'$)

$$\begin{aligned}\|\boldsymbol{\delta}_1^\pi(\mathbf{x}) - \mathbf{x}_1\| &\le \frac{\iint \|(\mathbf{I}+\mathbf{D})^{-1}\mathbf{z}_1\| g(\mathbf{c}, \mathbf{D}, \mathbf{z}_1, \mathbf{z}_2, \ldots, \mathbf{z}_m)\, d\mathbf{D}\, d\mathbf{H}}{\iint g(\mathbf{c}', \mathbf{D}, \mathbf{z}_1, \mathbf{z}_2, \ldots, \mathbf{z}_m)\, d\mathbf{D}\, d\mathbf{H}} \\ &\le C\frac{\iint \|(\mathbf{I}+\mathbf{D})^{-1}\mathbf{y}\| g(\mathbf{c}, \mathbf{D}, \mathbf{z}_1, \mathbf{z}_2, \ldots, \mathbf{z}_m)\, d\mathbf{D}\, d\mathbf{H}}{\iint g(\mathbf{c}', \mathbf{D}, \mathbf{z}_1, \mathbf{z}_2, \ldots, \mathbf{z}_m)\, d\mathbf{D}\, d\mathbf{H}},\end{aligned}$$

where $\mathbf{y}$ and $g$ are as in Lemma 3.8, with $u = a_1$ and $v = a_2 - a_1 + (m-1)/2$.

Now Lemma 3.8 shows that, if $u + v = a_2 + (m-1)/2 > 1$ and $u = a_1 < 1$, then

$$\int \|(\mathbf{I}+\mathbf{D})^{-1}\mathbf{y}\| g(\mathbf{c}, \mathbf{D}, \mathbf{z}_1, \mathbf{z}_2, \ldots, \mathbf{z}_m)\, d\mathbf{D} \le K \int g(\mathbf{c}', \mathbf{D}, \mathbf{z}_1, \mathbf{z}_2, \ldots, \mathbf{z}_m)\, d\mathbf{D}.$$

Hence $\|\boldsymbol{\delta}_1^\pi(\mathbf{x}) - \mathbf{x}_1\| \le K$ and the theorem is established. □

THEOREM 3.10. *Assume that $\pi(\boldsymbol{\beta})$ is $N_k(\boldsymbol{\beta}^0, \mathbf{A})$, $k \ge 2$, and $\pi(\mathbf{H}, \mathbf{D})$ satisfies Condition 1. Also suppose that we choose $l = 0$. If $a_1 < 1$ and $a_2 > 1 - \frac{m}{2}$, then $\boldsymbol{\delta}^\pi(\mathbf{x}) - \mathbf{x}$ is uniformly bounded.*

PROOF. The proof is very similar to that in Theorem 3.9:

$$\begin{aligned}\|(\nabla m(\mathbf{x}))_1\| \le \iint &\|(\mathbf{I}+\mathbf{V})^{-1}(\mathbf{x}_1 - \bar{\mathbf{x}}) + (\mathbf{I}+\mathbf{V}+m\mathbf{A})^{-1}\bar{\mathbf{x}}\| \\ &\times \frac{1}{|\mathbf{I}+\mathbf{D}|^{[a_2 - a_1 + (m-1)/2]}|\mathbf{I}+\mathbf{D}+m\mathbf{H}\mathbf{A}\mathbf{H}^t|^{1/2}|\mathbf{D}|^{a_1}} \\ &\times \exp\left(-\frac{1}{2}\sum_{i=1}^{m}(\mathbf{x}_i - \bar{\mathbf{x}})^t \mathbf{H}^t(\mathbf{I}+\mathbf{D})^{-1}\mathbf{H}(\mathbf{x}_i - \bar{\mathbf{x}})\right) \\ &\times \exp\left(-\frac{1}{2}m(\mathbf{H}\bar{\mathbf{x}})^t(\mathbf{I}+\mathbf{D}+m\mathbf{H}\mathbf{A}\mathbf{H}^t)^{-1}(\mathbf{H}\bar{\mathbf{x}})\right) d\mathbf{D}\, d\mathbf{H}.\end{aligned}$$

Note that

$$\begin{aligned}&\|(\mathbf{I}+\mathbf{V})^{-1}(\mathbf{x}_1 - \bar{\mathbf{x}}) + (\mathbf{I}+\mathbf{V}+m\mathbf{A})^{-1}\bar{\mathbf{x}}\| \\ &\quad \le \|(\mathbf{I}+\mathbf{V})^{-1}(\mathbf{x}_1 - \bar{\mathbf{x}})\| + \|(\mathbf{I}+\mathbf{V})^{-1}\bar{\mathbf{x}}\|.\end{aligned} \tag{38}$$

One now proceeds as in the proof of Theorem 3.9 with each term of (38), making use of Lemma 3.6 and arguments similar to the proof in that lemma. □

THEOREM 3.11. *Assume that* $\pi(\boldsymbol{\beta})$ *is* $N_k(\boldsymbol{\beta}^0, \lambda\mathbf{A})$, $k \geq 2$, $\pi(\mathbf{H}, \mathbf{D})$ *satisfies Condition* 1 *and* $\pi(\lambda)$ *satisfies Condition* 2. *Also suppose that we choose* $l = 0$. *If* $a_1 < 1$, $a_2 > 1 - \frac{m}{2}$ *and* $b > 1 - \frac{k}{2}$, *then* $\boldsymbol{\delta}^\pi(\mathbf{x}) - \mathbf{x}$ *is uniformly bounded.*

PROOF. Define $\boldsymbol{\delta}^\pi(\mathbf{x}|\lambda)$ to be the posterior mean with $\lambda$ given. From Theorem 3.10, we know that

$$\sup_{\mathbf{x}} \|\boldsymbol{\delta}^\pi(\mathbf{x}|\lambda) - \mathbf{x}\| \equiv K(\lambda) < \infty.$$

With a modification of the proof of Theorem 3.10, it can be shown that $K(\lambda)$ is continuous. Also, as $\lambda \to \infty$, the posterior distribution converges to that corresponding to $\pi(\boldsymbol{\beta}) = 1$, so we know from Theorem 3.9 that $\lim_{\lambda\to\infty} K(\lambda) < \infty$. As $\lambda \to 0$, the posterior converges to the special case of Theorem 3.10 in which $\mathbf{A} = \mathbf{0}$, so we know $K(0) < \infty$. It follows that $K(\lambda)$ is itself bounded. Finally, letting $\pi(\lambda|\mathbf{x})$ denote the posterior distribution of $\lambda$ given $\mathbf{x}$, which was shown to exist under the given conditions, it is clear that

$$\begin{aligned}\|\boldsymbol{\delta}^\pi(\mathbf{x}) - \mathbf{x}\|^2 &= \|E^{\pi(\lambda|\mathbf{x})}[\boldsymbol{\delta}^\pi(\mathbf{x}|\lambda) - \mathbf{x}]\|^2 \\ &\leq E^{\pi(\lambda|\mathbf{x})}\|\boldsymbol{\delta}^\pi(\mathbf{x}|\lambda) - \mathbf{x}\|^2 \leq E^{\pi(\lambda|\mathbf{x})}[K(\lambda)]^2.\end{aligned}$$

Since $K(\lambda)$ is bounded, it follows that $\|\boldsymbol{\delta}^\pi(\mathbf{x}) - \mathbf{x}\|$ is uniformly bounded, and the proof is complete. □

3.4. *Admissibility and inadmissibility results.* To prove admissibility or inadmissibility based on Results 3.1 and 3.2, we need only determine whether (25) is infinite or (26) is finite. Since Lemmas 3.5, 3.6 and 3.7 provide effectively equivalent upper and lower bounds on $m(\mathbf{x})$, it suffices to evaluate (25) and (26) for these equivalent bounds.

3.4.1. *Case* 1 *scenario.*

THEOREM 3.12. *Assume that* $\pi(\boldsymbol{\beta}) = 1$, $m \geq 2$, $a_1 < 1$ *and* $\pi(\mathbf{H}, \mathbf{D})$ *satisfies Condition* 1 *with* $l = 0$. *If* $k = 2$ *and* $a_2 > 1$, *then the posterior mean is admissible under quadratic loss. If* $\frac{3-m}{2} < a_2 < \frac{3}{2} - \frac{1}{k}$, *then the posterior mean is inadmissible.*

PROOF. Let $\mathbf{z}_i = (z_{i1}, z_{i2}, \ldots, z_{ik}) = \mathbf{H}(\mathbf{x}_i - \bar{\mathbf{x}})$. Define $y_j^2 = \sum_{i=1}^m z_{ij}^2$. By Lemma 3.5,

$$m(\mathbf{x}) \approx C \iint \frac{1}{|\mathbf{D}|^{a_1} |\mathbf{I} + \mathbf{D}|^{(m-1)/2 + a_2 - a_1}}$$

$$\times \exp\left(-\frac{1}{2}\sum_{i=1}^{m}(\mathbf{x}_i-\bar{\mathbf{x}})^t\mathbf{H}^t(\mathbf{I}+\mathbf{D})^{-1}\mathbf{H}(\mathbf{x}_i-\bar{\mathbf{x}})\right)\cdot d\mathbf{D}\,d\mathbf{H}$$

$$= C\iint\prod_{j=1}^{k}\frac{1}{d_j^{a_1}(1+d_j)^{(m-1)/2+a_2-a_1}}\exp\left(-\frac{y_j^2}{2(1+d_j)}\right)d\mathbf{D}\,d\mathbf{H}$$

$$= C\int\left[\prod_{j=1}^{k}\int_0^\infty\frac{1}{d_j^{a_1}(1+d_j)^{(m-1)/2+a_2-a_1}}\exp\left(-\frac{y_j^2}{2(1+d_j)}\right)dd_j\right]d\mathbf{H}.$$

Applying the upper bound of Lemma 3.3 with $r=(m-1)/2+a_2-a_1$ and $a=a_1$ yields

$$\int_0^\infty\frac{1}{d_j^{a_1}(1+d_j)^{(m-1)/2+a_2-a_1}}\exp\left(-\frac{y_j^2}{2(1+d_j)}\right)dd_j$$
$$\approx C_1\min\{C_2,(y_j^2)^{(3-m)/2-a_2}\}.$$

Thus

$$m(\mathbf{x})\approx C\int\left[\prod_{j=1}^{k}C_1\min\{C_2,(y_j^2)^{(3-m)/2-a_2}\}\right]d\mathbf{H}. \tag{39}$$

To prove admissibility, note that

$$\begin{aligned}\overline{m}(r) &= \int_{\{\mathbf{x}:\,\|\mathbf{x}\|=r\}} m(\mathbf{x})\,d\phi(\mathbf{x})\\ &\le C\int\left[\int_{\{\mathbf{x}:\,\|\mathbf{x}\|=r\}}\prod_{j=1}^{k}C_1\min\{C_2,(y_j^2)^{(3-m)/2-a_2}\}\,d\phi(\mathbf{x})\right]d\mathbf{H}.\end{aligned} \tag{40}$$

The inner integral, with respect to $\phi$, is essentially considering $\mathbf{x}$ to be uniformly distributed on the surface of the sphere of radius $\|\mathbf{x}\|=r$. Since $\mathbf{H}$ is an orthonormal matrix, $((\mathbf{H}\mathbf{x}_1)^t,(\mathbf{H}\mathbf{x}_2)^t,\ldots,(\mathbf{H}\mathbf{x}_m)^t)^t$ also has a uniform distribution on the surface of the sphere of radius $r$. From the result in Section 49, Subsection 1, of [24], it follows that, for each given $\mathbf{H}$,

$$\left(\frac{y_1^2}{r^2},\frac{y_2^2}{r^2},\ldots,\frac{y_k^2}{r^2},1-\frac{1}{r^2}\sum_{i=1}^{k}y_i^2\right)\sim\text{Dirichlet}\left(\frac{m-1}{2},\frac{m-1}{2},\ldots,\frac{m-1}{2},\frac{k}{2}\right).$$

Thus,

$$\overline{m}(r)\le C\iint_{\{\sum_{i=1}^{k}y_i^2\le r^2\}}\prod_{i=1}^{k}\min\{C_2,(y_i^2)^{(3-m)/2-a_2}\}$$
$$\times\prod_{i=1}^{k}\left(\frac{y_i^2}{r^2}\right)^{(m-3)/2}\left(1-\frac{1}{r^2}\sum_{i=1}^{k}y_i^2\right)^{(k-2)/2}d\left(\frac{y_1^2}{r^2}\right)\cdots d\left(\frac{y_k^2}{r^2}\right)d\mathbf{H}.$$

The inner integral is clearly constant over $\mathbf{H}$ and can be dropped, along with the factor $(1-\sum_{i=1}^{k} y_i^2/r^2)$ (since $k \geq 2$). Then elimination of the range restriction on the $y_i^2$ yields

$$\begin{aligned}\overline{m}(r) &\leq C \prod_{i=1}^{k} \int_0^{\infty} \min\{C_2, (y_i^2)^{[(3-m)/2-a_2]}\} \cdot \left(\frac{y_i^2}{r^2}\right)^{(m-5)/2} d\left(\frac{y_i^2}{r^2}\right) \\ &\leq C r^{-k(m-1)} \prod_{i=1}^{k} \int_0^{\infty} [C_2^{-1} + (y_i^2)^{-[(3-m)/2-a_2]}]^{-1} (y_i^2)^{(m-3)/2}\, dy_i^2,\end{aligned}$$

the last inequality using the fact that $\min\{C_2, v\} \leq 2(C_2^{-1} + v^{-1})^{-1}$. The final integrals are finite if $m \geq 2$ and $a_2 > 1$, so then

$$\overline{m}(r) \leq C r^{-k(m-1)}.$$

Hence

$$\int_c^{\infty} [r^{mk-1}\overline{m}(r)]^{-1}\, dr \geq \int_c^{\infty} \frac{1}{r^{k-1}}\, dr,$$

which is infinite if $k = 2$. Since the conditions $k = 2$ and $a_2 > 1$ also imply that $\boldsymbol{\delta}^{\pi}(\mathbf{x}) - \mathbf{x}$ is bounded by Theorem 3.9, the proof of admissibility using Result 3.1 is complete.

To prove inadmissibility, note from (39) that

$$\begin{aligned}\underline{m}(r) &= \int_{\{\mathbf{x}:\, \|\mathbf{x}\|=r\}} \frac{1}{m(\mathbf{x})}\, d\phi(\mathbf{x}) \\ &\leq \int_{\{\mathbf{x}:\, \|\mathbf{x}\|=r\}} \left(\int \prod_{j=1}^{k} C_1 \min\{C_2, (y_j^2)^{(3-m)/2-a_2}\}\, d\mathbf{H}\right)^{-1} d\phi(\mathbf{x}).\end{aligned} \tag{41}$$

Note that

$$\left[\int f(\mathbf{H}) d\mathbf{H}\right]^{-1} \leq \int [f(\mathbf{H})]^{-1}\, d\mathbf{H} \qquad \text{if } f(\mathbf{H}) > 0,$$

so that

$$\begin{aligned}\underline{m}(r) &\leq C_1 \int \left[\int_{\{\mathbf{x}:\, \|\mathbf{x}\|=r\}} \left(\prod_{j=1}^{k} \min\{C_2, (y_j^2)^{(3-m)/2-a_2}\}\right)^{-1} d\phi(\mathbf{x})\right] d\mathbf{H} \\ &\leq C_1 \int \left[\int_{\{\mathbf{x}:\, \|\mathbf{x}\|=r\}} \prod_{j=1}^{k} \max\{C_2, (y_j^2)^{a_2-(3-m)/2}\}\, d\phi(\mathbf{x})\right] d\mathbf{H}.\end{aligned} \tag{42}$$

Continuing exactly as in the proof of admissibility (but employing the bound $y_i^2 \le r^2$) yields

$$
\begin{aligned}
\underline{m}(r) &\le C \prod_{i=1}^{k} \int_0^{r^2} \max\{C_2, (y_i^2)^{a_2-(3-m)/2}\} \left(\frac{y_i^2}{r^2}\right)^{(m-3)/2} d\left(\frac{y_i^2}{r^2}\right) \\
&\le C \prod_{i=1}^{k} \int_0^{r^2} [C_2 + (y_i^2)^{a_2-(3-m)/2}] \left(\frac{y_i^2}{r^2}\right)^{(m-3)/2} d\left(\frac{y_i^2}{r^2}\right) \qquad (43) \\
&\le C + C_2 r^{k(2a_2+m-3)}.
\end{aligned}
$$

Hence

$$
\int_c^\infty r^{(1-mk)} \underline{m}(r)\, dr \le C + C_2 \int_c^\infty r^{(2ka_2-3k+1)}\, dr,
$$

which is finite only if $a_2 < \frac{3}{2} - \frac{1}{k}$. If $m \ge 2$, $a_1 < 1$ and $a_2 > (3-m)/2$, then $\boldsymbol{\delta}^\pi(\mathbf{x}) - \mathbf{x}$ is uniformly bounded, and Result 3.2 completes the proof of inadmissibility. (It was not strictly necessary to establish the uniform boundedness condition for inadmissibility, but it is necessary to verify that the posterior mean exists, and the uniform boundedness condition clearly establishes that this is so.) □

Theorem 3.12 fails to cover the situation in which $k = 2$ and $a_2 = 1$ and the situation $k \ge 3$ and $a_2 \ge \frac{3}{2} - \frac{1}{k}$. We suspect that the posterior mean is also inadmissible in these two situations, but were unable to prove it. (The main hurdle is to find a way to avoid use of the too-strong inequality $[\int f(\mathbf{H})\, d\mathbf{H}]^{-1} \le \int [f(\mathbf{H})]^{-1}\, d\mathbf{H}$.)

3.4.2. *Case 2 scenario.*

THEOREM 3.13. *Assume that* $\pi(\boldsymbol{\beta})$ *is* $N_k(\boldsymbol{\beta}^0, \mathbf{A})$, $a_1 < 1$, $k \ge 2$ *and* $\pi(\mathbf{H}, \mathbf{D})$ *satisfies Condition* 1 *with* $l = 0$. *If* $a_2 \ge 1 - \frac{1}{k}$, *then the posterior mean is admissible. If* $\frac{2-m}{2} < a_2 < 1 - \frac{1}{k}$, *then the posterior mean is inadmissible.*

PROOF. Let $\mathbf{z}_i = (z_{i1}, z_{i2}, \ldots, z_{ik}) = \mathbf{H}\mathbf{x}_i$. Define $y_j^2 = \sum_{i=1}^m z_{ij}^2$. By (33) we have

$$
m(\mathbf{x}) \approx C_1 \iint \frac{1}{|\mathbf{D}|^{a_1} |C_2\mathbf{I} + \mathbf{D}|^{m/2+a_2-a_1}} \times \exp\left(-\frac{1}{2}\sum_{i=1}^m \mathbf{x}_i^t \mathbf{H}^t (C_3\mathbf{I} + \mathbf{D})^{-1} \mathbf{H}\mathbf{x}_i\right) \cdot d\mathbf{D}\, d\mathbf{H}
$$

$$= C_1 \iint \prod_{j=1}^{k} \frac{1}{d_j^{a_1}(C_2+d_j)^{m/2+a_2-a_1}} \exp\left(-\frac{y_j^2}{2(C_3+d_j)}\right) d\mathbf{D}\, d\mathbf{H}$$

$$= C_1 \int \left[\prod_{j=1}^{k} \int_0^\infty \frac{1}{d_j^{a_1}(C_2+d_j)^{m/2+a_2-a_1}} \exp\left(-\frac{y_j^2}{2(C_3+d_j)}\right) dd_j\right] d\mathbf{H}.$$

Applying the upper bound of Lemma 3.3 with $r = a_2 - a_1 + m/2$ to the inner integral above yields

$$\int_0^\infty \frac{1}{d_j^{a_1}(C_2+d_j)^{m/2+a_2-a_1}} \exp\left(-\frac{y_j^2}{2(C_3+d_j)}\right) dd_j \le C_1^* \min\{C_2^*, (y_j^2)^{1-m/2-a_2}\}.$$

Thus

$$m(\mathbf{x}) \approx C \int \left[\prod_{j=1}^{k} C_1^* \min\{C_2^*, (y_j^2)^{1-m/2-a_2}\}\right] d\mathbf{H}. \tag{44}$$

To prove admissibility, note that

$$\begin{aligned} \overline{m}(r) &= \int_{\{\mathbf{x}:\, \|\mathbf{x}\|=r\}} m(\mathbf{x})\, d\phi(\mathbf{x}) \\ &\le C \int \left[\int_{\{\mathbf{x}:\, \|\mathbf{x}\|=r\}} \prod_{j=1}^{k} C_1^* \min\{C_2^*, (y_j^2)^{1-m/2-a_2}\}\, d\phi(\mathbf{x})\right] d\mathbf{H}. \end{aligned} \tag{45}$$

The inner integral with respect to $\phi$ is essentially considering $\mathbf{x}$ to be uniformly distributed on the surface of the sphere of radius $\|\mathbf{x}\| = r$. Since $\mathbf{H}$ is an orthonormal matrix, $((\mathbf{H}\mathbf{x}_1)^t, (\mathbf{H}\mathbf{x}_2)^t, \ldots, (\mathbf{H}\mathbf{x}_m)^t)^t$ also has a uniform distribution on the surface of the sphere of radius $r$. From the result in Section 49, Subsection 1, of [24], it follows that, for each given $\mathbf{H}$,

$$\left(\frac{y_1^2}{r^2}, \frac{y_2^2}{r^2}, \ldots, \frac{y_k^2}{r^2}\right) \sim \text{Dirichlet}\left(\frac{m}{2}, \frac{m}{2}, \ldots, \frac{m}{2}\right).$$

Thus,

$$\overline{m}(r) \le C \iint \prod_{i=1}^{k} \min\{C_2^*, (y_i^2)^{1-m/2-a_2}\} \cdot \prod_{i=1}^{k} \left(\frac{y_i^2}{r^2}\right)^{m/2-1} d\left(\frac{y_i^2}{r^2}\right) d\mathbf{H}.$$

Again dropping the integral over $\mathbf{H}$ and using the inequality $\min\{C_2, v\} \le 2(C_2^{-1} + v^{-1})^{-1}$ results in the bound

$$\begin{aligned} \overline{m}(r) &\le C \prod_{i=1}^{k} \int_0^{r^2} [C_2^{*-1} + (y_i^2)^{-(1-m/2-a_2)}]^{-1} \left(\frac{y_i^2}{r^2}\right)^{m/2-1} d\left(\frac{y_i^2}{r^2}\right) \\ &\le C r^{-km} \left[\int_0^{r^2} (C_2^{*-1} + v^{(a_2+m/2-1)})^{-1} v^{(m-2)/2}\, dv\right]^k. \end{aligned}$$

The order of the integral in the last expression is easily seen to be $O(r^{2(1-a_2)})$ if $a_2 < 1$; $O(\log r)$ if $a_2 = 1$; and $O(1)$ if $a_2 > 1$. Hence

$$\int_c^\infty [r^{mk-1}\overline{m}(r)]^{-1}\,dr \ge \begin{cases} C\int_c^\infty r^{1-2k(1-a_2)}\,dr, & \text{if } a_2 < 1,\\ C\int_c^\infty r(\log r)^{-k}\,dr, & \text{if } a_2 = 1,\\ C\int_c^\infty r\,dr, & \text{if } a_2 > 1.\end{cases}$$

This is clearly infinite if $a_2 \ge 1 - 1/k$. By Theorem 3.10, this condition also implies that $\boldsymbol{\delta}^\pi(\mathbf{x}) - \mathbf{x}$ is bounded, so use of Result 3.1 completes the proof of admissibility.

To prove inadmissibility, note from (44) and the fact $[\int f(\mathbf{H})\,d\mathbf{H}]^{-1} \le \int [f(\mathbf{H})]^{-1}\,d\mathbf{H}$ that

$$\begin{aligned}\underline{m}(r) &= \int_{\{\mathbf{x}\colon \|\mathbf{x}\|=r\}} \frac{1}{m(\mathbf{x})}\,d\phi(\mathbf{x})\\ &\le \int_{\{\mathbf{x}\colon \|\mathbf{x}\|=r\}} \left(\int \prod_{j=1}^k C_1 \min\{C_2, (y_j^2)^{1-m/2-a_2}\}\,d\mathbf{H}\right)^{-1} d\phi(\mathbf{x})\\ &\le C\int\left[\int_{\{\mathbf{x}\colon \|\mathbf{x}\|=r\}} \prod_{j=1}^k \max\{C_2, (y_j^2)^{a_2+(m-2)/2}\}\,d\phi(\mathbf{x})\right]d\mathbf{H}.\end{aligned}$$

Continuing as with the proof of admissibility, one obtains

$$\begin{aligned}\underline{m}(r) &\le C\prod_{i=1}^k \int_0^{r^2} \max\{C_2, (y_i^2)^{a_2-(m-2)/2}\}\left(\frac{y_i^2}{r^2}\right)^{(m-2)/2} d\left(\frac{y_i^2}{r^2}\right)\\ &\le C\prod_{i=1}^k \int_0^{r^2} [C_2 + (y_i^2)^{a_2-(m-2)/2}]\left(\frac{y_i^2}{r^2}\right)^{(m-2)/2} d\left(\frac{y_i^2}{r^2}\right)\\ &\le C + C_2 r^{k(2a_2+m-2)}.\end{aligned}$$

Hence

$$\int_c^\infty r^{(1-mk)}\underline{m}(r)\,dr \le C + C_2\int_c^\infty r^{(2ka_2-2k+1)}\,dr,$$

which is finite only if $a_2 < 1 - 1/k$. If $a_1 < 1$ and $a_2 > (2-m)/2$, then $\boldsymbol{\delta}^\pi(\mathbf{x}) - \mathbf{x}$ is uniformly bounded and so the posterior mean exists, and Result 3.2 completes the proof of inadmissibility. $\square$

3.4.3. *Case 3 scenario.*

THEOREM 3.14. *Assume that* $\pi(\boldsymbol{\beta})$ *is* $N_k(\boldsymbol{\beta}^0, \lambda\mathbf{A})$, $m \geq 2$, $a_1 < 1$, $\pi(\mathbf{H}, \mathbf{D})$ *satisfies Condition* 1 *with* $l = 0$ *and* $\pi(\lambda)$ *satisfies Condition* 2. *If* (i) $k \geq 2$, $a_2 \geq 1 - \frac{1}{k}$ *and* $b > 1$; *or* (ii) $k \geq 3$, $a_2 > 1 - \frac{b}{k}$ *and* $0 \leq b < 1$; *or* (iii) $k = 2$, $a_2 > 1 - \frac{b}{2}$ *and* $0 < b < 1$, *then the posterior mean is admissible under quadratic loss.*

PROOF. Starting with (35) from Lemma 3.7 (setting all constants to 1 for notational simplicity) yields

$$m(\mathbf{x}) \approx \iiint \prod_{j=1}^{k} \frac{1}{d_j^{a_1}(1+d_j)^{(m-1)/2+a_2-a_1}(1+\lambda+d_j)^{1/2}} \times \exp\left(-\frac{1}{2}\left[\sum_{i=1}^{m} \frac{(\mathbf{H}(\mathbf{x}_i - \bar{\mathbf{x}}))_j^2}{1+d_j} + m\frac{(\mathbf{H}\bar{\mathbf{x}})_j^2}{1+\lambda+d_j}\right]\right) \times \pi(\lambda)\, d\lambda\, d\mathbf{D}\, d\mathbf{H}.$$

Let

$$w_j = \sum_{i=1}^{m} \frac{(\mathbf{H}(\mathbf{x}_i - \bar{\mathbf{x}}))_j^2}{\|\mathbf{x}\|^2}, \qquad v_j = m\frac{(\mathbf{H}\bar{\mathbf{x}})_j^2}{\|\mathbf{x}\|^2}, \qquad j = 1, 2, \ldots, k.$$

Under $\phi(\mathbf{x})$, the uniform distribution on the surface of the sphere of radius $r = \|\mathbf{x}\|$, by the result in Section 49, Subsection 1, of [24], we have

$$(w_1, \ldots, w_k, v_1, \ldots, v_k) \sim \mathrm{Dirichlet}\left(\frac{m-1}{2}, \ldots, \frac{m-1}{2}, \frac{1}{2}, \ldots, \frac{1}{2}\right).$$

Thus, arguing as in previous theorems, dropping $\mathbf{H}$ and letting the $w_i$ and $v_i$ range freely over $(0, 1)$, yields

$$\begin{aligned}\overline{m}(r) &= \int_{\|\mathbf{x}\|=r} m(\mathbf{x})\, d\phi(\mathbf{x}) \\ &\leq \int \prod_{j=1}^{k} \left\{ \int_0^1 \int_0^1 \int_0^1 \frac{1}{d_j^{a_1}(1+d_j)^{(m-1)/2+a_2-a_1}(1+\lambda+d_j)^{1/2}} \times \exp\left(-\frac{r^2}{2}\left[\frac{w_j}{1+d_j} + \frac{v_j}{1+\lambda+d_j}\right]\right) \times w_j^{(m-1)/2-1} v_j^{-1/2}\, dw_j\, dv_j\, dd_j \right\} \pi(\lambda)\, d\lambda.\end{aligned}$$

Make the change of variables $s_j = w_j(1+d_j)$ and $t_j = v_j(1+d_j+\lambda)$, $j = 1, 2, \ldots, k$. The region of integration becomes

$$R_{st} = \left\{0 < s_j \leq \frac{1}{1+d_j}, 0 < t_j \leq \frac{1}{1+d_j+\lambda}, i = 1, 2, \ldots, k\right\}.$$

Then

$$\begin{aligned}\overline{m}(r) \le \iiiint_{R_{st}} \prod_{j=1}^{k} \biggl\{ &\frac{((1+d_j)s_j)^{(m-3)/2}((1+d_j+\lambda)t_j)^{-1/2}}{d_j^{a_1}(1+d_j)^{(m-1)/2+a_2-a_1}(1+\lambda+d_j)^{1/2}} \\ &\times \exp\biggl(-\frac{r^2}{2}[s_j+t_j]\biggr) \\ &\times (1+d_j)(1+d_j+\lambda)\,ds_j\,dt_j\,dd_j\biggr\}\pi(\lambda)\,d\lambda \\ \le \iiiint_{R_{st}} \prod_{j=1}^{k} \biggl\{ &\frac{1}{d_j^{a_1}(1+d_j)^{a_2-a_1}} s_j^{(m-3)/2} t_j^{-1/2} \\ &\times \exp\biggl(-\frac{r^2}{2}[s_j+t_j]\biggr)\,ds_j\,dt_j\,dd_j\biggr\}\pi(\lambda)\,d\lambda \\ = \int \prod_{j=1}^{k} \biggl\{ \int &\frac{1}{d_j^{a_1}(1+d_j)^{a_2-a_1}} \\ &\times \biggl[\int_0^{1/(1+d_j)} s_j^{(m-3)/2}\exp\biggl(-\frac{r^2}{2}s_j\biggr)\,ds_j\biggr] \\ &\times \biggl[\int_0^{1/(1+d_j+\lambda)} t_j^{-1/2}\exp\biggl(-\frac{r^2}{2}t_j\biggr)\,dt_j\biggr]\,dd_j\biggr\}\pi(\lambda)\,d\lambda.\end{aligned}$$

Applying Lemma 3.3(b) to the inner integrals above yields

$$\begin{aligned}\overline{m}(r) \le C\int \prod_{j=1}^{k} \biggl[\int &\frac{1}{d_j^{a_1}(1+d_j)^{a_2-a_1}} \\ &\times \min\{r^{-(m-1)}, (1+d_j)^{-(m-1)/2}\} \\ &\times \min\{r^{-1}, (1+d_j+\lambda)^{-1/2}\}\,dd_j\biggr]\pi(\lambda)\,d\lambda.\end{aligned} \tag{46}$$

Consider first the situation $b>1$. Then $\pi(\lambda)$ has finite mass and so [using $(1+d_j+\lambda)^{-1/2} \le (1+d_j)^{-1/2}$]

$$\begin{aligned}\overline{m}(r) \le C\biggl[\int_0^\infty &\frac{1}{d^{a_1}(1+d)^{a_2-a_1}}\min\{r^{-(m-1)}, (1+d)^{-(m-1)/2}\} \\ &\times \min\{r^{-1}, (1+d)^{-1/2}\}\,dd_j\biggr]^k.\end{aligned}$$

Break up the inner integral into integrals $I_1$ and $I_2$ over $(0, r^2-1)$ and $(r^2-1, \infty)$, respectively. Then, since $a_1 < 1$,

$$I_1 = \int_0^{r^2-1} \frac{1}{d^{a_1}(1+d)^{a_2-a_1}} \cdot \frac{1}{r^{(m-1)}} \cdot \frac{1}{r}\,dd \le Cr^{-m}(1+r^{2(1-a_2)}),$$

$$I_2 = \int_{r^2-1}^{\infty} \frac{1}{d^{a_1}(1+d)^{a_2-a_1}} \cdot \frac{1}{(1+d)^{(m-1)/2}} \cdot \frac{1}{(1+d)^{1/2}}\, dd \le Cr^{2-2a_2-m}.$$

Hence

$$\overline{m}(r) \le C[I_1 + I_2]^k \le Cr^{-mk}(1 + r^{2k(1-a_2)})$$

and

$$\int_c^{\infty} [r^{mk-1}\overline{m}(r)]^{-1}\, dr \ge \int_c^{\infty} \frac{r}{1 + r^{2k(1-a_2)}}\, dr.$$

This is finite if $a_2 \ge 1 - 1/k$.

Next consider the case $0 \le b \le 1$. Clearly

$$\begin{aligned} &\min\{r^{-1}, (1 + d_j + \lambda)^{-1/2}\} \\ &\quad = \min\{r^{-1}, (1 + d_j + \lambda)^{-1/2}\}^{2(1-b)+\varepsilon} \\ &\qquad \times \min\{r^{-1}, (1 + d_j + \lambda)^{-1/2}\}^{2b-1-\varepsilon} \\ &\quad \le (1+\lambda)^{(b-1-\varepsilon/2)} \min\{r^{-1}, (1 + d_j)^{-1/2}\}^{2b-1-\varepsilon}. \end{aligned}$$

Hence (46) can be bounded as

$$\begin{aligned} \overline{m}(r) \le C\biggl[\int_0^{\infty} &\frac{1}{d^{a_1}(1+d)^{a_2-a_1}} \\ &\times \min\{r^{-(m-1)}, (1+d)^{-(m-1)/2}\} \\ &\qquad \times \min\{r^{-1}, (1+d)^{-1/2}\}\, dd\biggr]^{(k-1)} \\ \times \int_0^{\infty} &\frac{1}{d^{a_1}(1+d)^{a_2-a_1}} \\ &\times \min\{r^{-(m-1)}, (1+d)^{-(m-1)/2}\} \\ &\times \min\{r^{-1}, (1+d)^{-1/2}\}^{(2b-1-\varepsilon)}\, dd \end{aligned}$$

[using the fact that $(1+\lambda)^{(b-1-\varepsilon/2)}\pi(\lambda)$ has finite mass]. Proceeding exactly as in the $b > 1$ case yields

$$\overline{m}(r) \le Cr^{[2k-2ka_2-km+2(1-b)+\varepsilon]},$$

so that

$$\int_c^{\infty} [r^{mk-1}\overline{m}(r)]^{-1}\, dr \ge C\int_c^{\infty} r^{(2k-2ka_2+1-2b+\varepsilon)}\, dr,$$

which is infinite if $a_2 \ge 1 - \frac{b}{k} + \varepsilon'$. Since $\varepsilon'$ was arbitrary, the condition for admissibility when $0 \le b < 1$ is $a_2 > 1 - \frac{b}{k}$. By Theorem 3.11 these conditions also imply that $\boldsymbol{\delta}^{\pi}(\mathbf{x}) - \mathbf{x}$ is uniformly bounded, except when $k = 2$, in which case the restriction $b > 0$ must be added. This completes the proof of admissibility. □

3.4.4. *Admissibility and inadmissibility for the common priors.* Let us apply these results to the versions of the reference prior discussed in the Introduction. For $\boldsymbol{\beta}$, the Case 1 constant prior leads to admissibility only in the case $k=2$, and hence is not a prior we recommend. The Case 2 conjugate prior can readily yield admissible estimators, and is certainly reasonable if backed by subjective knowledge. The Case 3 default prior that was suggested in Section 1.2 is

$$\pi(\boldsymbol{\beta}) \propto [1 + \|\boldsymbol{\beta}\|^2]^{-(k-1)/2}, \tag{47}$$

corresponding to the two-stage prior $\boldsymbol{\beta}|\lambda \sim N(\mathbf{0}, \lambda\mathbf{I})$, $\pi(\lambda) \propto \lambda^{-1/2} e^{-1/(2\lambda)}$. We therefore focus on admissibility results when this prior is used for $\boldsymbol{\beta}$.

In regard to priors for $\mathbf{V}$, note first that the nonhierarchical reference prior for $\mathbf{V}$ cannot be considered, since it corresponds to $a_1 = 1$, yielding an improper posterior. The modification

$$\pi(\mathbf{V}) = \frac{1}{|\mathbf{V}|^{a_1} \prod_{i<j}(d_i - d_j)},$$

where $a_1 < 1$, is inadmissible in Case 1 [$\pi(\boldsymbol{\beta}) = 1$], but is admissible in Case 2, and is admissible in Case 3 when $b > 1$ and $1 - 1/k \le a_1 < 1$, or when $0 < b < 1$ and $1 - b/k < a_1 < 1$. Since we recommend (47), which has $b = 1/2$, this suggests the choice $a_1 = 1 - 1/(2k) = (2k-1)/(2k)$. While we were not strictly able to prove admissibility for this choice, it likely corresponds to admissibility and, in any case, being at the boundary of admissibility has considerable appeal.

The modified reference prior of the form

$$\pi(\mathbf{V}) = \frac{1}{|\mathbf{I} + \mathbf{V}|^{a_2} \prod_{i<j}(d_i - d_j)}$$

is admissible in Case 1 if $k = 2$ and $a_2 > 1$; in Case 2 or Case 3 ($b > 1$) if $a_2 \ge 1 - 1/k$; and in Case 3 ($0 < b < 1$) if $a_2 > 1 - b/k$. The natural choice is $a_2 = 1$, since this is admissible for all $b$ and $k$ in Cases 2 and 3, and is almost admissible in Case 1 when $k = 2$. (Recall that we were unable to establish admissibility or inadmissibility in this case, but again being at the boundary of admissibility has considerable appeal.) Recalling the discussion from the Introduction, the recommended default prior distribution of this form is thus

$$\pi(\mathbf{V}) = \frac{1}{|\mathbf{I} + \mathbf{V}| \prod_{i<j}(d_i - d_j)}.$$

This yields a proper posterior with a posterior mean that is admissible in estimation under quadratic loss.

## APPENDIX

PROOF OF LEMMA 3.3. (a) It suffices to take $c_1 = c_2 = 1$ in the proof. This is because $(c_1 + d)/(c_2 + d)$ is uniformly bounded above and below, so that one could change $(c_1 + d)$ to $(c_2 + d)$, or vice versa. A simple change of variables then reduces the expression to the case $c_2 = 1$. Clearly,

$$\begin{aligned} f(v) &= \int_0^1 \frac{1}{(1+d)^r d^a} \exp\left(-\frac{v}{2(1+d)}\right) dd \\ &\quad + \int_1^\infty \frac{1}{(1+d)^r d^a} \exp\left(-\frac{v}{2(1+d)}\right) dd \\ &\le e^{-v/4} \int_0^1 \frac{1}{d^a}\, dd + \int_1^\infty \frac{1}{d^{r+a}} \exp\left(-\frac{v}{4d}\right) dd \\ &= \frac{1}{1-a} e^{-v/4} + \int_1^\infty \frac{1}{d^{r+a}} \exp\left(-\frac{v}{4d}\right) dd. \end{aligned}$$

Making the change of the variables $t = v/d$ yields

$$\begin{aligned} f(v) &= \frac{1}{1-a} e^{-v/4} + \int_v^0 \left(\frac{t}{v}\right)^{r+a} \exp\left(-\frac{t}{4}\right)\left(-\frac{v}{t^2}\right) dt \\ &= \frac{1}{1-a} e^{-v/4} + v^{1-r-a} \int_0^v t^{r+a-2} e^{-t/4}\, dt \\ &\le \frac{1}{1-a} e^{-v/4} + v^{1-r-a} \int_0^\infty t^{(r+a-1)-1} e^{-t/4}\, dt \\ &= \frac{1}{1-a} e^{-v/4} + v^{1-r-a}\Gamma(r+a-1) \cdot 4^{r+a-1}. \end{aligned}$$

Since $r + a > 1$, it is easy to show that $e^{-v/4} \le C v^{1-r-a}$ when $v \ge 0$. Therefore

$$f(v) \le C_1 v^{1-r-a}.$$

On the other hand, $f(v)$ is a decreasing function of $v$ when $v \ge 0$, so

$$\begin{aligned} \max f(v) = f(0) &= \int_0^\infty \frac{1}{(1+d)^r d^a}\, dd \\ &\le \int_0^1 \frac{1}{d^a}\, dd + \int_1^\infty \frac{1}{d^{r+a}}\, dd = \frac{1}{1-a} + \frac{1}{r+a-1} = C_2. \end{aligned}$$

Thus, defining $C_3 = C_2/C_1$,

$$f(v) \le C_1 \min\{C_3, v^{1-r-a}\}.$$

To find a lower bound for $f(v)$, note that

$$\begin{aligned} f(v) &\geq \int_1^\infty \frac{1}{(1+d)^r d^a} \exp\left(-\frac{v}{2(1+d)}\right) dd \\ &\geq \int_1^\infty \frac{1}{(2d)^r d^a} \exp\left(-\frac{v}{2d}\right) dd \\ &= \frac{1}{2^r} \int_1^\infty \frac{1}{d^{r+a}} \exp\left(-\frac{v}{2d}\right) dd. \end{aligned}$$

Making the change of the variables $t = v/d$, one obtains

$$\begin{aligned} f(v) &\geq \frac{1}{2^r} \int_v^0 \left(\frac{t}{v}\right)^{r+a} \exp\left(-\frac{t}{2}\right)\left(-\frac{v}{t^2}\right) dt \\ &= \frac{1}{2^r} v^{1-r-a} \int_0^v t^{r+a-2} \exp\left(-\frac{t}{2}\right) dt. \end{aligned}$$

If $v \geq 1$, then

$$f(v) \geq \frac{1}{2^r} v^{1-r-a} \int_0^1 t^{r+a-2} \exp\left(-\frac{t}{2}\right) dt = C_1' v^{1-r-a}.$$

If $0 \leq v < 1$, then

$$f(v) \geq f(1) = \int_0^\infty \frac{1}{(1+d)^r d^a} \exp\left(-\frac{1}{2(1+d)}\right) dd = C_2'.$$

Let $C_3' = \min\{C_1', C_2'\}$. Thus

$$f(v) \geq \min\{C_3', C_1' v^{1-r-a}\} = C_1' \min\{C_4', v^{1-r-a}\},$$

where $C_4' = C_3'/C_1'$, completing the proof of part (a).

To prove part (b), change variables from $t$ to $w = xt$. Then

$$g(\mu, v) = v^{-(a+1)} \int_0^{\mu v} w^a e^{-w} \, dw.$$

Now $\int_0^{\mu v} w^a e^{-w} \, dw < \Gamma(a+1)$ and $\int_0^{\mu v} w^a e^{-w} \, dw < \int_0^{\mu v} w^a \, dw = (\mu v)^{(a+1)}/(a+1)$. Hence

$$g(\mu, v) < v^{-(a+1)} \min\left\{\Gamma(a+1), \frac{(\mu v)^{(a+1)}}{(a+1)}\right\},$$

and the result follows. □

J. O. Berger
ISDS
Box 90251
Duke University
Durham, North Carolina 27708-0251
USA
e-mail: berger@stat.duke.edu

W. E. Strawderman
Department of Statistics
Rutgers University
Hill Center, Busch Campus
110 Frelinghuysen Avenue
Piscataway, New Jersey 08854-8019
USA
e-mail: straw@stat.rutgers.edu

D. Tang
Biostatistics and Statistical Report
Novartis Pharmaceuticals
One Health Plaza
East Hanover, New Jersey 07936
USA
e-mail: dejun.tang@pharma.novartis.com